%% file: main.tex
\documentclass[11pt]{article}

\input{preambles}

\begin{document}

\def\paperTitle{Fast and Accurate Intersections on a Sphere}
\headers{\paperTitle}{Hongyu Chen, Paul A. Ullrich, and Julian Panetta}

\title{\paperTitle}

\author{Hongyu Chen\thanks{Computer Science Graduate Group, University of California, Davis, Davis, CA 
  (\email{hyvchen@ucdavis.edu}).}
\and Paul A. Ullrich\thanks{Physical and Life Science Division, Lawrence Livermore National Laboratory, Livermore, CA and Department of Land, Air and Water Resources, University of California, Davis, Davis, CA
  (\email{paullrich@ucdavis.edu}).}
\and Julian Panetta\thanks{Department of Computer Science,  University of California, Davis, Davis, CA
  (\email{jpanetta@ucdavis.edu}).}
}
\date{} 
\maketitle

\let\thefootnote\relax
\footnotetext{\textbf{Funding:} This work was funded by the U.S. Department of Energy Office of Science Program for Climate Model Diagnosis and Intercomparison (PCMDI) at Lawrence Livermore National Laboratory under master award DE-AC52-07NA27344.}

\input{paper_body}  

\clearpage
\phantomsection
\addcontentsline{toc}{section}{Supplementary Material}
\thispagestyle{plain}
\begin{center}
  {\Large\bfseries \MakeUppercase{Supplementary Material}}
\end{center}
\vspace{1ex}
\markboth{SUPPLEMENTARY MATERIAL}{SUPPLEMENTARY MATERIAL}

\makeatletter
\let\ACCUX@orig@section\section
\renewcommand{\section}{\subsection}
\makeatother

\setcounter{subsection}{0}
\setcounter{equation}{0}\setcounter{figure}{0}\setcounter{table}{0}
\numberwithin{equation}{subsection}
\numberwithin{figure}{subsection}
\numberwithin{table}{subsection}
\renewcommand{\thesubsection}{S\arabic{subsection}}
\renewcommand{\theequation}{S\arabic{subsection}.\arabic{equation}}
\renewcommand{\thefigure}{S\arabic{subsection}.\arabic{figure}}
\renewcommand{\thetable}{S\arabic{subsection}.\arabic{table}}
\crefname{subsection}{Supplement}{Supplements}

\input{supplement_body}  

\bibliographystyle{plainnat} 
\bibliography{references}

\end{document}

%% file: preambles.tex
\usepackage[letterpaper,margin=0.9in,includehead,includefoot]{geometry}
\usepackage[labelfont=bf,font=it]{caption}

\usepackage[utf8]{inputenc} 
\usepackage[T1]{fontenc}
\usepackage{lmodern}
\usepackage{microtype}

\usepackage{amsmath, amsfonts, amssymb}
\usepackage{mathtools}
\usepackage{amsthm}
\usepackage{graphicx}
\usepackage{epstopdf} 
\usepackage{booktabs}
\usepackage{algorithm}
\usepackage{algorithmic}
\usepackage{xcolor}
\usepackage{xspace}
\usepackage[numbers,sort&compress]{natbib}

\usepackage[bookmarks=false]{hyperref} 
\usepackage[nameinlink]{cleveref}

\crefname{equation}{Eq.}{Eqs.}
\crefname{figure}{Figure}{Figures}
\crefname{table}{Table}{Tables}
\crefname{section}{Section}{Sections}
\crefname{subsection}{Section}{Sections}
\crefname{algorithm}{Algorithm}{Algorithms}

\newenvironment{keywords}{%
  \par\smallskip\noindent\textbf{Keywords.}\ }{\par\smallskip}

\newenvironment{MSCcodes}{%
  \par\smallskip\noindent\textbf{MSC Codes.}\ }{\par\smallskip}

\newcommand{\email}[1]{\href{mailto:#1}{#1}}

\newcommand{\headers}[2]{}

\numberwithin{equation}{section}
\theoremstyle{plain}

\newcommand{\newsiamthm}[2]{\newtheorem{#1}{#2}[section]}
\theoremstyle{remark}
\newcommand{\newsiamremark}[2]{\newtheorem{#1}{#2}}
\newsiamremark{remark}{Remark}
\newsiamremark{hypothesis}{Hypothesis}
\crefname{hypothesis}{Hypothesis}{Hypotheses}
\newsiamthm{claim}{Claim}

\usepackage[normalem]{ulem}

\renewcommand{\vec}[1]{{\bf #1}}
\providecommand{\cross}{\times}
\def\n{\vec{n}}
\def\nhat{\hat{\vec{n}}}
\def\x{\vec{x}}
\def\p{\vec{p}}
\def\e{\vec{e}}

\def\t{\vec{t}}
\def\P{\vec{P}}

\newcommand{\defeq}{\vcentcolon=}
\newcommand{\eqdef}{=\vcentcolon}
\providecommand{\norm}[1]{{\left\|#1\right\|}}
\providecommand{\abs}[1]{\lvert#1\rvert}

\providecommand{\UB}[1]{{\left[#1\right]}_\text{UB}}
\providecommand{\LB}[1]{{\left[#1\right]}_\text{LB}}

\makeatletter
\newcommand{\fltimes}{\mathbin{\mathpalette\make@circled *}}
\newcommand{\flplus}{\mathbin{\mathpalette\make@circled {\raisebox{.2\height}{\scalebox{.7}{+}}} }}
\newcommand{\flsub}{\mathbin{\mathpalette\make@circled {\raisebox{.2\height}{\scalebox{.7}{\(-\)}}} }}
\newcommand{\fldiv}{\mathbin{\mathpalette\make@circled {\raisebox{.2\height}{\scalebox{.7}{\( \hspace{0.3em} \divisionsymbol \) }}} }}
\newcommand{\flstar}{\mathbin{\mathpalette\make@circled {\raisebox{.3\height}{\scalebox{.7}{\( \hspace{0.325em} \star \) }}} }}
\newcommand{\make@circled}[2]{%
  \ooalign{$\m@th#1\smallbigcirc{#1}$\cr\hidewidth$\m@th#1#2$\hidewidth\cr}%
}
\newcommand{\smallbigcirc}[1]{%
  \vcenter{\hbox{\scalebox{0.77778}{$\m@th#1\bigcirc$}}}%
}
\makeatother

\providecommand{\subexpression}{\mathcal{S}}
\providecommand{\flvalue}[1]{v_{#1}}
\providecommand{\flerror}[1]{\varepsilon_{#1}}
\providecommand{\flabserror}[1]{\Delta_{#1}}

\providecommand{\hot}{\text{h.o.t.}}
\let\divisionsymbol\div

\def\ie{\emph{i.e.}\xspace}

\renewcommand{\matrix}[1]{\begin{bmatrix} #1 \end{bmatrix}}

%% file: paper_body.tex
\maketitle
\begin{abstract}
We introduce a fast, high-precision algorithm for calculating intersections
between great circle arcs and lines of constant latitude on the unit sphere. We
first propose a simplified intersection point formula with improved speed and
numerical robustness over the ones traditionally implemented in geoscience
software. We then show how algorithms based on the concept of error-free
transformations (EFT) can be applied to evaluate this formula within a relative
error bound that is on the order of machine precision.
We demonstrate that, with a vectorized and parallelized implementation, this
enhanced accuracy is achieved with no compute-time overhead compared to a
direct calculation in hardware floating point, making our algorithm suitable for
performance-sensitive applications like regridding of high-resolution climate
data. In contrast, evaluating our formula using high-precision data types like
quadruple precision and arbitrary precision, or using the robust intersection
computation routines from the Computational Geometry Algorithms Library (CGAL),
leads to significant computational overhead, especially since these alternatives
inhibit vectorization. More generally, our work demonstrates how EFT techniques
can be combined and extended to implement nontrivial geometric calculations with
high accuracy and speed.
\end{abstract}



\section{Introduction}
Calculating points of intersection between arcs defined on a sphere is an
important subproblem of meshing algorithms used in geoscience. For instance,
regridding algorithms \cite{TempestRemap,strategiesRemapping,FARRELL20092632} compute the intersection of two different
meshes of the globe.
The two most common variants of the arc intersection problem involve
intersections between (1) two geodesic arcs, or (2) a geodesic arc and a line of
constant latitude. Geodesic arcs trace along the circle formed by intersecting
the sphere with a plane passing through the origin, while lines of constant
latitude are defined by intersecting the sphere with the plane $z = z_0$. Note that
every line of constant \emph{longitude} is a geodesic arc, but the equator is
the only line of constant latitude that is also a geodesic.

The first type of intersection is relatively simple: it can be computed by
evaluating three cross products and normalizing the result. Fast and accurate
cross-product implementations based on Kahan's $2 \times 2$ determinant
algorithm \cite{kahancost} are well known, and tight bounds on their errors have
been established \cite{KahanAlgo}; these implementations have already
been adopted by geoscience software packages like UXarray \cite{uxarray2024} and
Climate Data Operators (CDO) \cite{CDO} to reduce roundoff errors
in geodesic intersection calculations.
The second type, however, is more complex. It has drawn our focus due the lack
of implementations offering both speed and accuracy guarantees and its
importance to common data analysis operations on the globe like computing zonal
means \cite{TempestRemap,TempestRemap2,strategiesRemapping}. While deriving
closed-form expressions for intersection points is straightforward
(\Cref{sec:math}), evaluating these accurately in finite precision and bounding
the resulting error is nontrivial.

For the high-resolution inputs typical in geoscience applications, performance
is critical, leading the software packages developed in those communities like
UXarray and CDO to calculate intersections of the second type in the hardware's
native floating-point arithmetic. However, numerical stability issues that arise
when evaluating the associated formulas cause these packages to produce results
that are imprecise---sometimes dramatically so as seen in our experiments.
At the other extreme, computational geometry software like the 3D
Spherical Geometry Kernel \cite{cgal:cclt-sgk3-24b,cgal:eb-24b} of CGAL (the
Computational Geometry Algorithms Library) \cite{cgal:eb-24b} can produce
\emph{exact} intersection coordinates represented in algebraic number types. As
shown in our timing comparisons, this approach is several orders of magnitude
slower than plain floating-point calculation, making it unsuitable for
large-scale computations. Furthermore, downstream algorithms in geoscience
computations generally operate on floating-point representations of the
intersection coordinates, requiring the exactly computed results to be replaced
by rounded approximations in the end.

We show that, by evaluating an improved intersection point formula using
an approach that extends and combines error-free transformation (EFT) techniques
\cite{Ogita2005}, it is possible to achieve accuracy approaching that of the
double-precision approximation computed by CGAL with no compute time
overhead versus the imprecise hardware-arithmetic code (in a batch-processing
setting; see \Cref{sec:experiments}). While our algorithm
executes more arithmetic operations than a literal translation of the formula
into floating point instructions, both computations reach the limits of the hardware's memory
bandwidth after they are fully parallelized using multithreading and
vectorization. This high level of performance is critical for the intended
application of large-scale meshing and data processing algorithms.

We establish the accuracy of our algorithm by proving that, to first order, the
relative error in the calculated point is bounded by a small multiple of the
machine epsilon, $u$, that is independent of the input data
provided a mild nondegeneracy condition holds. Without this condition,
the leading term remains $O(u)$, but with a larger input-dependent constant.
This is in contrast to a literal floating-point implementation whose error we
show is of order $O(\sqrt{u})$ for challenging inputs, indicating a loss of
roughly half the working precision that is observed in our experiments.
We also demonstrate our EFT-based implementation's
accuracy on a large suite of randomly generated test cases, comparing against
evaluating the same intersection point formulas using two standard
high-precision arithmetic libraries, GNU MPFR \cite{MPFR} (implementing
arbitrary-precision floating-point types) and GCC's \texttt{libquadmath}
(quadruple-precision IEEE 754 arithmetic) \cite{libquadmath}.
Our algorithm computes results of comparable accuracy to both libraries in
significantly less time, and this performance gap widens when vectorization is
enabled since the libraries' number types inhibit the use of SIMD instructions.


\subsection{Previous Work}
Detecting intersections between curves and calculating intersection point
coordinates are core problems of computational geometry. Despite the simplicity
of the involved formulas, robust implementations in inexact computer arithmetic
are nontrivial even in the simple case of lines in the plane. Significant effort
has been devoted to developing algorithms that offer both speed and robustness
guarantees.
A prominent example is the Triangle library \cite{Triangle96}, which uses
adaptive precision \cite{richard_shewchuk_adaptive_1997}, refining its
calculations as needed to ensure predicates like orientation and incircle tests
are always evaluated correctly. Exact predicate implementations like
these---coupled with perturbation strategies to break ties consistently
\cite{Edelsbrunner1990,Devillers2003}---typically suffice to ensure a geometric
algorithm successfully terminates and produces a valid output (\emph{e.g.}, a
true Delaunay triangulation). However, ensuring the accuracy of the geometric
quantities constructed by these algorithms is a separate issue. For example,
the Triangle library evaluates the coordinates of line segment intersections
inexactly.

Guaranteeing exactness of the numerical results computed by an algorithm
requires performing operations in a suitable number type. For instance, the
intersection of two line segments whose endpoint coordinates are floating-point
numbers is a point with rational coordinates.
The CGAL library implements exact computations when the user selects one of its
kernels with ``exact constructions,'' which for certain algorithms and inputs is
needed to avoid crashes. Most relevant to our problem, CGAL's 3D Spherical
Geometry Kernel \cite{cgal:cclt-sgk3-24b,cgal:eb-24b} computes exact
intersections between arbitrary circles using a special data type representing
algebraic numbers of degree two. However, when floating point representations
are needed downstream, these exact results must be approximated, and we show in
\Cref{fig:comparison_methods,fig:pole_area_accuracy} that the approximation
routines provided by CGAL do not provide significant accuracy improvements over
our much faster evaluations.

Developing robust and practical software for computing intersections (and
other Boolean operations) between high-resolution triangle and tetrahedral
meshes in 3D remains an active area of research in computer graphics
\cite{Zhou:2016:MASG,Hu2018TetMeshWild,hu_fast_2020,Cherchi2020,trettner_ember_2022}.
With the exception of \texttt{fTetWild} \cite{hu_fast_2020},
these recent methods perform all of their constructions exactly (or implicitly)
and require a final approximation step to obtain a floating-point output.
This conversion raises accuracy concerns as noted above
and can introduce inverted or degenerated elements in difficult cases. The
\texttt{fTetWild} algorithm \cite{hu_fast_2020}, in contrast, does all
constructions in floating point, perturbing all output vertices as needed to
ensure predicates are satisfied; this approach is fast and guarantees a valid
floating-point output, though it sacrifices accuracy of the intersection point
computations. None of these works addresses meshing problems on the sphere, but
they serve as inspiration for future algorithms that we hope to develop, for
which the accurate floating point calculations studied here are a key
ingredient. CGAL does support computing arrangements of great circle arcs on the
sphere and other parametric surfaces
\cite{cgal:wfzh-a2-24b,berberich2010arrangements}, but the arrangements package
currently does not support inserting lines of constant latitude. Furthermore,
the result would require careful post-processing to obtain accurate
floating-point results.

Existing software developed in the geoscience community like UXarray
\cite{uxarray2024}, CDO \cite{CDO}, and CDO's underlying YAC (Yet Another Coupler) library \cite{YAC}
eschew exact number types in favor of
floating-point arithmetic for performance reasons.
These codes adopt some
strategies for reducing roundoff errors, like the accurate cross-product
implementation mentioned above, but the accuracy of many of their computations has not
been carefully investigated.
We show that, for the specific intersection point
calculation problem studied here, substantial gains in accuracy can be achieved
over these implementations with no sacrifice of performance.
A recent report by Mirshak \cite{mirshak2024intersections} presents a
solution to our intersection point problem based on Krumm's method for
intersecting planes \cite{krumm2000};
the Python-based floating point implementation it proposes suffers from
severe numerical issues as we show in \Cref{sec:experiments}.

The simplest way to increase the precision of a floating-point
calculation is to switch from the double-precision arithmetic natively supported
in hardware to a high-precision number type. The GNU MPFR library \cite{MPFR} is
a widely used implementation of arbitrary-precision floating-point computation,
enabling users to define new types with specified numbers of precision bits.
In fact, MPFR and its underlying GNU Multi-Precision Arithmetic Library
(GMP) \cite{GMP} are both used internally by CGAL.
Each arithmetic operation performed with MPFR types is correctly rounded
to the specified precision following the rules from the IEEE 754 standard \cite{ieee754},
enabling standard error analysis techniques to be applied to the results.
Another more limited but more performant option is to use the IEEE 754
quadruple-precision floating-point type (binary128). While this generally
lacks hardware support, a software implementation is available in GCC via the
\texttt{libquadmath} library \cite{libquadmath}. Our work shows that both of
these approaches introduce significant performance overhead and inhibit
vectorization.

An orthogonal strategy for improving computational accuracy is employing the
Newton-Raphson method to iteratively refine the calculated intersection
point from an initial floating point approximation. However, implementing this method with
standard floating-point arithmetic does not automatically guarantee improved accuracy
\cite{DennisWalker1984}. Additionally, running multiple Newton iterations
can impact performance.


To improve accuracy without major compromises to performance, we leverage Error-Free
Transformations (EFT) \cite{rump2009error} to simulate higher
precision where needed. The basic principle of this approach is that the
roundoff error committed by a primitive floating-point operation like addition
or multiplication is itself a floating-point number that can be computed with
simple algorithms \cite{Knuth1997,Graillat2009,Ogita2005}. Therefore, the exact
result of the operation can be represented as the sum of two floating
point numbers, where
the first holds the rounded value and the second holds the roundoff error. This
pair can then be propagated through subsequent operations to ensure that the
roundoff error is compensated for. This idea underlies the adaptive precision
scheme \cite{richard_shewchuk_adaptive_1997} used to evaluate predicates in the
Triangle library, which represents exact values as the sum of one or more
non-overlapping floating point numbers, called components,
ordered by magnitude. A benefit of this number format for predicate evaluation
is that it supports lazy evaluation: computation can begin on the highest-order
component and terminate once the result is sufficiently accurate to determine
the predicate's value. Our work employs just a two-component scheme, which
suffices to obtain highly accurate results across our large-scale test suite.
EFT-based algorithms have been proposed to accurately
evaluate dot products \cite{Ogita2005}, sums of squares \cite{Graillat2015}, and
square roots \cite{Rump2023} with rigorous error bounds, and we use these algorithms as
subroutines in our intersection point calculation.

\begin{table}
\footnotesize
\caption{Notation}\label{tab:notation}
    \vspace{-0em}
\begin{center}
  \renewcommand{\arraystretch}{1.2}
  \begin{tabular}{|c|l|} \hline
   \bf Notation & \bf Description \\ \hline
    $\delta_i$                    & Relative error incurred by floating-point operators ($+,-,*,\divisionsymbol$), with $|\delta_i| \le u$. \\ \hline
    $\flabserror{\subexpression}$ & Absolute error in the value of subexpression $\subexpression$ as computed by a particular \\ & algorithm being analyzed: $\flabserror{\subexpression} = \flvalue{\subexpression} - \subexpression$, where $\flvalue{\subexpression}$ is the computed value. \\ \hline
    $\flerror{\subexpression}$    & Relative error in the value of subexpression/quantity $\subexpression$ as computed by a \\ & particular algorithm being analyzed: $\flerror{\subexpression} = \flabserror{\subexpression} / S. $ \\ \hline
      $\tilde \subexpression = \bar{\subexpression} + e_{\bar{\subexpression}}$
                                  & Enhanced-accuracy result computed for quantity $\subexpression$, represented as the \\
                                  & unevaluated sum of $(\bar{\subexpression}, e_{\bar{\subexpression}})$, with leading part $\bar{\subexpression}$ and error compensation $e_{\bar{\subexpression}}$. \\ \hline
    $\mathbf{n}=(n_x, n_y, n_z)$  & Exact cross product \( \mathbf{n} = \mathbf{x}_1 \times \mathbf{x}_2 \). \\ \hline
  \end{tabular}
\end{center}
\end{table}

\subsection{Problem Definition}
\begin{figure}[H]
    \centering
    \includegraphics[width=0.6\linewidth]{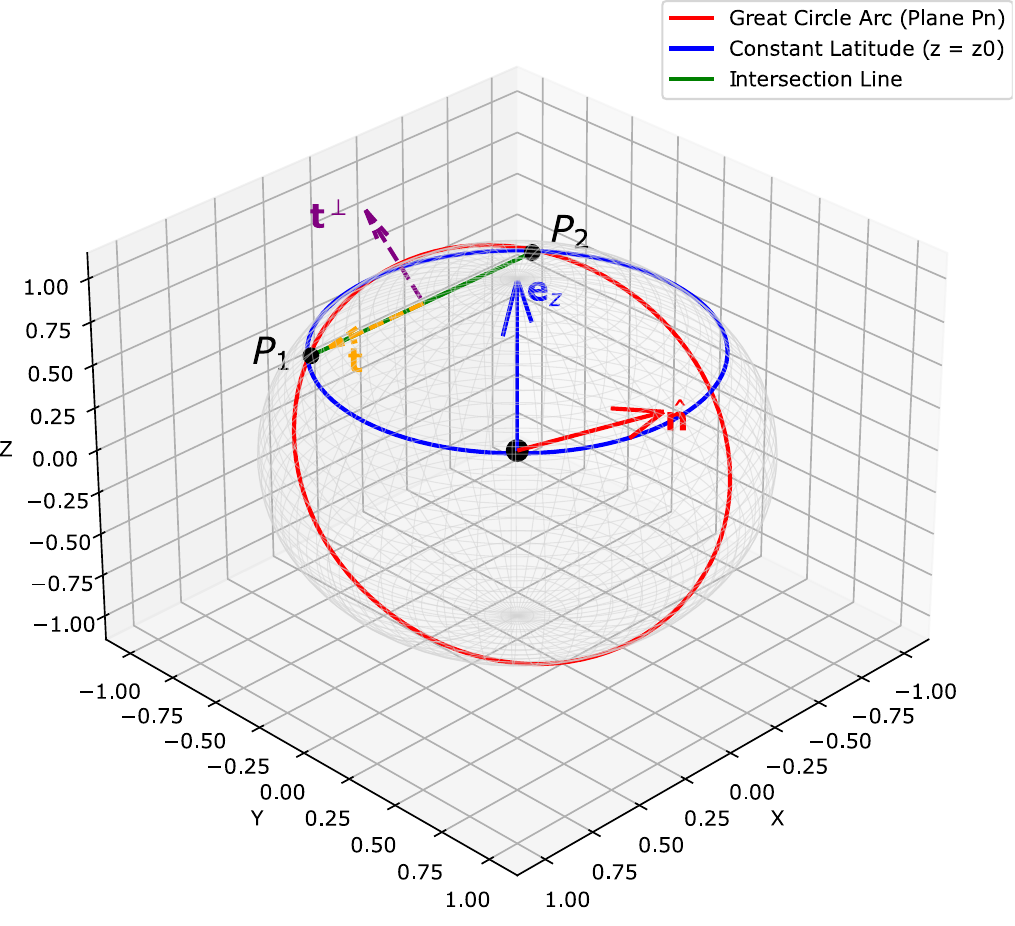} 
    \caption{
        The red circle represents the great circle arc defined by the plane \( P_{\hat{n}} \), which passes through the origin with normal vector \( \hat{n} \).
        The constant-latitude circle is drawn in blue.
        The vector \( \hat{e}_z \) is aligned with the z-axis.
        Also visualized are the vectors \( \t \) and \( \t^\perp \) introduced in \Cref{sec:math}, which are tangent and perpendicular to the planes' intersection line, respectively.
    }
    \label{fig:unit_sphere_visualization}
\end{figure}

A geodesic arc, or great circle arc, is the shortest path between two points \((\mathbf{x}_1, \mathbf{x}_2)\) on the unit sphere. It can be constructed by intersecting the sphere with the plane passing through the origin with normal vector
\(
\mathbf{n} = \mathbf{x}_1 \times \mathbf{x}_2
\).
A constant-latitude circle on the sphere is defined by intersecting with the plane \(z = z_0\), representing a horizontal slice at \(z_0\).
Depending on the parameters of the arcs, they may intersect at 0, 1, or 2 points.
Given inputs $\x_1$, $\x_2$, and $z_0$ in floating point\footnote{
Because the coordinates of inputs $\x_1$ and $\x_2$ are represented by
approximate floating point numbers, these points are unlikely to lie precisely on the unit sphere.
However, since our intersection formula
\cref{eq:final_gca_constLat_intersection_point} is invariant to $\norm{\x_1}$
and $\norm{\x_2}$,
our algorithm can be interpreted as
operating on sphere points $\frac{\x_1}{\norm{\x_1}}$ and $\frac{\x_2}{\norm{\x_2}}$ despite
those generally being unrepresentable in floating point.
}, we seek to compute floating-point
approximations to these intersection point(s) as illustrated in \Cref{fig:unit_sphere_visualization}.
Determining analytical formulas for the intersection point coordinates is
straightforward (\Cref{sec:math}), however we will show that the na\"ive
translation of these formulas into floating-point instructions done in
existing geoscience software can lead to inaccurate results.

To assess the accuracy with which an intersection point
is computed, we measure the Euclidean distance between it
and the corresponding exact intersection \(\P = [P_x, P_y, z_0]\),
computing a relative error:
\begin{equation}
\label{eq:error_assessment}
    \frac{\norm{\flabserror{\P}}}{\norm{\P}} = \norm{\flabserror{\P}} = \sqrt{{(\flabserror{P_x})}^2 + {(\flabserror{P_y})}^2},
\end{equation}
where the denominator vanishes because \(\P\) lies on the unit sphere.
Note that while we establish bounds for the relative error of the intersection point as a whole,
the individual coordinate relative errors $\flerror{P_x}$ and $\flerror{P_x}$ are undefined if \(P_x = 0\) or \(P_y\ = 0\), respectively.

We develop a method for computing intersection points that achieves
high relative accuracy while incurring no performance penalty
compared to existing implementations when batch-processing large datasets.
We show that our method, when
implemented in a floating point arithmetic with \( p \) precision bits, computes solutions
with accuracy
\[
\norm{\flabserror{\P}} \leq \varphi u + O(u^2),
\]
where \(\varphi\) is a moderate constant that is independent of the input provided a mild nondegeneracy assumption holds,
and \( u = 2^{-p} \) is the unit roundoff. This bound means that
our method computes results of comparable accuracy to
those obtained by employing high-precision data types (or exact computation)
and rounding coordinates back to precision $p$ at the end.

\subsection{Outline of the Paper}

The paper is organized as follows. In \Cref{sec:math},
we derive a standard analytical intersection point formula
and present a simplified version that is faster to evaluate and easier to analyze.
In \Cref{sec:notation}, we provide some background material on
floating point arithmetic and introduce our notation. \Cref{sec:float_error} presents a numerical analysis of our proposed math formula's direct floating point implementation and indicates the potential for significant loss of precision. 
\Cref{sec:accurate_operators} summarizes the EFT concept
and reviews existing algorithms for computing basic quantities
like dot products, and squared norms and difference of products.
\Cref{sec:accurate_intersection_point_calculation} details our
\texttt{AccuX} algorithm building on these techniques
and analyzes its accuracy.
\Cref{sec:experiments} documents the results of a large-scale experiment
we conducted to compare the accuracy and performance of our algorithm
to existing methods. Finally, \Cref{sec:conclusions} and \Cref{sec:future_work}
draw conclusions and identify directions for future work.

\section{Intersection Point Formulas}
\label{sec:math}
Here we derive a standard analytical formula for
the desired intersection point(s) that is implemented in \cite{YAC}.
We then identify a simple improvement to this formula that accelerates its evaluation
and simplifies its numerical analysis.
Intersections lie simultaneously on three surfaces:
the unit sphere, the $z = z_0$ plane, and the plane defining the geodesic arc
with unit normal $\nhat \defeq \frac{\x_1 \cross \x_2}{\norm{\x_1 \cross \x_2}} =
[\hat{n}_x, \hat{n}_y, \hat{n}_z]$.
The two planes intersect in a line perpendicular to their normals, which can
be expressed in parametric form as:
$$
\p(\alpha) = \alpha \t + z_0 \frac{\t^\perp}{\e_z \cdot \t^\perp},
$$
where $\t \defeq \nhat \cross \e_z$ is a vector tangent to the intersection line, with  $\e_z$  being the unit axis vector associated with the Cartesian \( z \)-coordinate. The vector $\t^\perp \defeq \t \cross \nhat$ represents a 90-degree rotation of $\t$  around $\nhat$.

This line intersects the sphere at parameters $\alpha^*$ solving the equation:
$$
1 = \norm{\p(\alpha)}^2 = {\alpha}^2 \norm{\t}^2 + z_0^2 \frac{\norm{\t^\perp}^2}{(\e_z \cdot \t^\perp)^2}
                          = {\alpha}^2 \norm{\t}^2 + z_0^2 \frac{1}{\norm{\t}^2},
$$
where we simplified $\norm{\t^\perp} = \norm{\t}$
and $\e_z \cdot \t^\perp = \e_z \cdot (\t \cross \nhat) = \t \cdot (\nhat \cross \e_z) = \norm{\t}^2$.
We note that $\norm{\t}^2 = \hat{n}_x^2 + \hat{n}_y^2$ and conclude:
$$
\alpha^* = \mp \frac{1}{\norm{\t}^2} \sqrt{\norm{\t}^2 - z_0^2}
\implies
\mathbf{P} = \p(\alpha^*) = \frac{1}{\hat{n}_x^2 + \hat{n}_y^2} \left( z_0 \t^\perp \mp \sqrt{\hat{n}_x^2 + \hat{n}_y^2 - z_0^2} \t \right).
$$

\noindent
Evaluating the cross product expressions for $\t$ and $\t^\perp$ obtains the coordinates:

\begin{equation}
\mathbf{P} = 
\begin{bmatrix}
- \frac{1}{\hat{n}_x^2 + \hat{n}_y^2}\left(z_0 \hat{n}_x \hat{n}_z \pm \hat{s} \hat{n}_y\right) \\
- \frac{1}{\hat{n}_x^2 + \hat{n}_y^2}\left(z_0 \hat{n}_y \hat{n}_z \mp \hat{s} \hat{n}_x\right) \\
z_0
\end{bmatrix},
\label{eq:gca_constlat_intersection_point_cdo}
\end{equation}
where \( \hat{s} \defeq \sqrt{\hat{n}_x^2 + \hat{n}_y^2 - z_0^2} \).
While we have not found this formula reported in the literature, it has been implemented in \cite{YAC}.

We observe that the normalization of $\nhat$ is unnecessary since both the numerator
and the denominator of each coordinate expression is proportional to its squared norm.
Thus \cref{eq:gca_constlat_intersection_point_cdo} can be rewritten using the non-unit
normal $\n \defeq \x_1 \cross \x_2 = [n_x, n_y, n_z]$ and its
prefix $\n_{xy} \defeq [n_x, n_y]$:

\begin{equation}
\mathbf{P} = 
\begin{bmatrix}
    - \frac{1}{\norm{\n_{xy}}^2}\left(z_0 n_x n_z \pm s n_y\right) \\
    - \frac{1}{\norm{\n_{xy}}^2}\left(z_0 n_y n_z \mp s n_x\right) \\
z_0
\end{bmatrix},
\label{eq:final_gca_constLat_intersection_point}
\end{equation}
with \( s = \sqrt{\norm{\n_{xy}}^2 - \|\mathbf{n}\|^2 z_0^2} \).
This simplified formula is roughly
1.5 times faster in our experiments (\Cref{sec:experiments})
and simplifies the error analysis since the roundoff error committed by
normalization no longer needs to be propagated through the full calculation.

The three surfaces intersect if and only if
\[
\norm{\n_{xy}}^2 \ge \|\mathbf{n}\|^2 z_0^2,
\]
and there exist two points of intersection when this inequality is strict.
We note that any intersections found must still be checked
to ensure they lie between $\x_1$ and $\x_2$.
In the following, we assume that intersections exist and focus on computing
their coordinates accurately. From a numerical accuracy perspective, the
two solutions are equivalent, so we specifically focus on computing
the $\alpha^* \le 0$ solution.

\section{Notation and Assumptions}
\label{sec:notation}
Our work assumes all floating-point operations follow the IEEE 754 standard \cite{ieee754}. Let $\mathbb{F}$ denote the set of floating-point numbers with base $\beta$ and precision $p$ (both integers):
\[
\mathbb{F} = \{0\} \cup \{M \cdot \beta^e : M, e \in \mathbb{Z}, \ \beta^{p-1} \le |M| < \beta^p \},
\]
where the exponent range is left unbounded to avoid considering underflow and overflow.
Our work specifically assumes binary floating point ($\beta = 2$), for which the unit roundoff error is $u = 2^{-p}$ ($2^{-53}$ for double precision).

We define $\operatorname{fl} : \mathbb{R} \rightarrow \mathbb{F}$ as a round-to-nearest function
\[
|t - \operatorname{fl}(t)| = \min_{f \in \mathbb{F}} |t - f|, \quad t \in \mathbb{R},
\]
without assuming a specific tie-breaking rule. Our computational experiments
employ IEEE 754's default rounding mode (round-to-nearest with ties broken by rounding to even $M$).
The IEEE 754 standard guarantees that the value computed in floating point for an
operation $a \circ b$ on floating-point inputs $a, b \in \mathbb{F}$
is $\operatorname{fl}(a \circ b)$ for any ${\circ \in \{+, -, *, /\}}$.
To be precise about the rounding steps in our algorithms,
we exclusively use $\operatorname{fl}(\cdot)$ to denote a single rounding;
in contrast to \cite{Higham2002} and \cite{Ogita2005}, we do not
use the overloaded notation $\operatorname{fl}(\subexpression)$
to refer to an evaluation of subexpression $\subexpression$ involving multiple rounding steps
(whose result can depend on the order of operations and whether
floating point contractions are employed).
Thus, in our notation, the expression $\operatorname{fl}(a * b - c)$ on
inputs $a, b, c \in \mathbb{F}$ implies
executing a fused multiply-subtraction operation, which commits only a single rounding.

%
Throughout this work, we assume that all inputs are floating point numbers:
\( {\mathbf{x}_1 = [x_1, y_1, z_1] \in \mathbb{F}^3} \),
\( {\mathbf{x}_2 = [x_2, y_2, z_2] \in \mathbb{F}^3}\),
and \(z_0 \in \mathbb{F}\).
We further consider the input geodesic arc to satisfy \(n_x, n_y, n_z \geq 0\) and the constant latitude to be in the northern hemisphere, \emph{i.e.}, \(z_0 \geq 0\); this assumption does not reduce generality because we can transform any geodesic arc with negative normal vector components and any line of constant latitude in the southern hemisphere by reflecting across the $x$, $y$, or $z$ planes as necessary.
\Cref{tab:notation} summarizes the notation we employ in our error analyses.

\section{Evaluating \Cref{eq:final_gca_constLat_intersection_point} in Floating Point}
\label{sec:float_error}
The most straightforward approach to computing intersection points is to
directly translate \Cref{eq:final_gca_constLat_intersection_point} into
floating-point operations. We provide an error
analysis for this algorithm in \Cref{sec:direct_fp_analysis} of the
supplementary material
%
Our analysis identifies different asymptotic error behavior
in two different regimes depending on the size of subexpression
$s^2 \defeq \norm{{\n_{xy}}}^2 - z_0^2 \norm{\n}^2$.

When $s^2 \approx |\flabserror{s^2}|$ (referred to as regime (i)),
our analysis yields an error bound of
\begin{equation}
\label{eq:native_float_regime_i}
\norm{\flabserror{\P}} \le
        \sqrt{\frac{21}{\norm{\n_{xy}}}} \sqrt{u} + \hot
\end{equation}
This indicates a potential loss of approximately half the working precision,
which is confirmed by our experiments (\Cref{fig:comparison_methods} and
\Cref{fig:comparison_methods_extreme}). The issue underlying this precision loss
is the square-root evaluation of $\sqrt{s^2}$, which ends up
amplifying $|\flabserror{s^2}|$ from $O(u)$ to $O(\sqrt{u})$ in this regime (and also when $s^2 < |\flabserror{s^2}|$).

When $s^2 \gg |\flabserror{s^2}|$ (referred to as regime (ii)),
we assume a certain level of magnitude separation,
$s^2 > |\flabserror{s^2}|^{1/c}$ for some $1 < c < \infty$, to obtain
the bound:
$$
    \norm{\flabserror{\P}} 
             \le \frac{21}{\sqrt{\norm{\n_{xy}}}^{\frac{1}{c}}} u^\frac{2c - 1}{2c} + \hot
$$
This indicates that precision is lowered by at worst a factor $\frac{2 c - 1}{2 c}$, with
no lowering in the limit of infinite separation $c \to \infty$.

Alternatively, we can derive a bound without these regime assumptions:
$$
\norm{\flabserror{\P}} \le \left(\frac{21}{s} + \frac{19 \sqrt{2}}{\norm{\n_{xy}}}\right) u + \hot,
$$
obtaining a leading-order term $u$ that blows up for inputs with $s \to 0$.

\section{Accurate Operators}
\label{sec:accurate_operators}
The $O(\sqrt{u})$ error bound obtained in \Cref{sec:float_error} suggests that achieving $O(u)$ accuracy for the intersection point requires simulating twice the precision of the underlying floating point type.
We do this using EFT-based operations. From the mathematical formulation in \Cref{sec:math}, we see that the key operations needed include cross products; the sum of squares for norm calculations; the dot product for computing generic sums of products; and the square root. 

EFT-based algorithms have been previously published for the dot product, sum of squares, and square root operations, which we adapt with minor modifications to support chaining these operations together to maintain accuracy across larger computations. To compute entries of the cross product vector, we propose an approach with greater accuracy than Kahan's widely used difference-of-products algorithm and derive an error bound for it. These developments, along with detailed explanations of the operators, are presented in the following subsections.

\subsection{Accurate Dot Product}

The \texttt{CompDot} algorithm (first introduced as \texttt{Dot2} in \cite{Ogita2005}) computes accurate dot products by compensating for rounding errors during intermediate steps. It uses two EFT building blocks: \texttt{TwoSum} \cite{Knuth1997,MullerEtAl2018} and \texttt{TwoProd} \cite{dekker1971floating,Ogita2005,LefevreEtAl2023,Rump2023} in its internal calculations and then returns its result by rounding to a single floating-point number.

To retain the higher precision with which this result was computed, we propose the composable algorithm, \texttt{CompDotC}, which outputs both the sum and the error compensation term as a pair. This avoids the rounding step, improving accuracy for downstream operations. This variant of the algorithm is presented in \Cref{alg:refined_dot_product}, while the original \texttt{CompDot} algorithm is listed in \cref{alg:compensated_dot_product}.

\begin{algorithm}
\caption{The proposed \texttt{CompDotC}}
\label{alg:refined_dot_product}
\begin{algorithmic}[1]
\STATE \textbf{Input:} $\vec{x}, \vec{y} \in \mathbb{F}^n$
\STATE Initialize $[p, s] \gets \text{TwoProduct}(\vec{x}[0], \vec{y}[0])$
\FOR{$i = 1$ to $n-1$}
    \STATE $[h, r] \gets \text{TwoProduct}(\vec{x}[i], \vec{y}[i])$
    \STATE $[p, q] \gets \text{TwoSum}(p, h)$
    \STATE $s_i \gets \text{fl}(s + (q + r))$
\ENDFOR
\STATE \textbf{Return:} $(p_n, s_n)$
\end{algorithmic}
\end{algorithm}

\begin{algorithm}[H]
\caption{Original Compensated Dot Product \texttt{CompDot} (Algorithm 5.4 in \cite{Ogita2005})}
\label{alg:compensated_dot_product}
\begin{algorithmic}[1]
\STATE \textbf{Input:} $\vec{x}, \vec{y} \in \mathbb{F}^n$
    \STATE $[p_n, s_n] \gets \text{CompDotC}(\vec{x}, \vec{y})$
\STATE $\text{res} \gets \text{fl}(p_n + s_n)$ 
\STATE \textbf{Return} $\text{res}$
\end{algorithmic}
\end{algorithm}
The following relative error bound for \texttt{CompDotC} is proved by Ogita et al. \cite{Ogita2005}:
\begin{equation}
    \label{eq:dotfma_re_rel_err}
    |\varepsilon_\texttt{CompDotC}| =
    \frac{| (\vec{x}^\top \vec{y} - p_n) - s_n |}{|\vec{x}^\top \vec{y}|}
    \leq \gamma_n \frac{n u}{1 - (n - 1)u}
    \frac{|\vec{x}| \cdot |\vec{y}|}{|\vec{x}^\top \vec{y}|}
    \eqdef
    \UB{\varepsilon_\texttt{CompDotC}},
\end{equation}
where
$\gamma_n := \frac{n u}{1 - n u}$ for $n \in \mathbb{N}$,
with the implicit assumption that $n u < 1$ \cite{Higham2002}. 

The relative error bound for \texttt{CompDot} includes an additional $u$ term accounting for the rounding error incurred by explicitly evaluating the pair sum:
\begin{align}
    \label{eq:dotfma_rel_err}
    \UB{\varepsilon_\texttt{CompDot}} \vcentcolon= u + \UB{\varepsilon_\texttt{CompDotC}}.
\end{align}

\subsection{Accurate Cross Product}
The cross product of vectors \(\mathbf{u}\) and \(\mathbf{v}\) is:

\[
\mathbf{u} \times \mathbf{v} = (u_y v_z - u_z v_y, u_z v_x - u_x v_z, u_x v_y - u_y v_x).
\]

\noindent The operation used to compute each component is a $2 \times 2$ determinant, also referred to as a difference of products (\texttt{DiffOfProd}). This operation suffers from catastrophic cancellation when \( u_y v_z \approx u_z v_y \), making accurate evaluation challenging.

A classical technique for evaluating $2 \times 2$ determinants with improved accuracy is Kahan's algorithm (Algorithm \ref{alg:Kahan}), which uses Fused Multiply-Add (FMA) operations to estimate and compensate for rounding errors, ensuring a relative error bound of \(2u\)
\cite{KahanAlgo}.

\begin{algorithm}[H]
\caption{Kahan's Algorithm for \texttt{DiffOfProd} \cite{KahanLectureNote}}
\label{alg:Kahan}
\begin{algorithmic}[1]
\STATE \textbf{Input:} $a, b, c, d \in \mathbb{F}$
\STATE $bc \gets \text{fl}(b * c)$
\STATE $\text{err} \gets \texttt{FMA}(-b, c, bc)$ 
\STATE $\text{dop} \gets \texttt{FMA}(a, d, -bc)$
\STATE \textbf{Return:} $\text{dop} + \text{err}$
\end{algorithmic}
\end{algorithm}

For our problem, \texttt{DiffOfProd} is an intermediate step, and ultimately we need its error to be of size $O(u^2)$ to achieve $O(u)$ accuracy in the calculated intersection point. To this end, we propose \texttt{AccuDOP} (\Cref{alg:refined_dot_product}), which we derived from the accurate dot product algorithm. It returns both the main result and its associated error compensation term as a pair rather than summing them together to avoid a final $O(u)$ roundoff error.

\begin{algorithm}[H]
\caption{Proposed \texttt{AccuDOP} \cite{KahanLectureNote, KahanAlgo}}
\label{alg:AccuDOP}
\begin{algorithmic}[1]
\STATE \textbf{Input:} $a, b, c, d \in \mathbb{F}$
    \STATE $[p_1, r_1] \gets \text{TwoProduct}(a, d)$
    \STATE $[p_2, r_2] \gets \text{TwoProduct}(b, c)$
    \STATE $ [\text{dop},s]\gets \text{TwoSum}(p_1, -p_2)$
    \STATE $\text{err} \gets \text{fl}(s + \text{fl}(r_1 + r_2))$
\STATE \textbf{Return:} $[\text{dop},\text{err}]$
\end{algorithmic}
\end{algorithm}
Using the existing bound from \Cref{alg:refined_dot_product} \cite{Ogita2005}, it is trivial to derive the following relative error bound for \texttt{AccuDOP}:
\[
\left| \frac{(\text{dop} + \text{err}) - (ad - bc)}{ad - bc} \right| \leq \left(1 + 2 \frac{\left(|ad| + |bc|\right)}{|ad - bc|}\right) \cdot u^2 + O(u^3)
\]
For cross products of unit vectors, we can simplify this bound by replacing \(|ad| + |bc|\) in the numerator with 1. For instance, consider the evaluation of $n_x$, whose magnitude is clearly bounded:
$$
 |y_1 z_2 - y_2 z_1| = |n_x| \le \norm{\n} \le \norm{\x_1} \norm{\x_2} \le 1.
$$
When $y_1 z_2$ and $y_2 z_1$ have different signs, this implies:
$$
|y_1 z_2| + |y_2 z_1| = |y_1 z_2 - y_2 z_1| \le 1.
$$
When the signs are the same, the first equality above does not hold, but instead:
$$
|y_1 z_2| + |y_2 z_1| = |y_1 z_2 + y_2 z_1| \le 1.
$$
The inequality above follows from interpreting $y_1 z_2 + y_2 z_1$ as the first component of the cross product of $\x_1$ with the vector obtained by reflecting $\x_2$ across either the $y=0$ or $z = 0$ plane. Since this is still the cross product of two unit vectors, its component magnitudes must not exceed 1.

The argument above can be repeated for $n_y$ and $n_z$, leading to a simplified relative error bound for each normal component:
\begin{equation}
\label{eq:accu_normal_error}
\left| \frac{(\bar{n}_{x/y/z} + e_{\bar{n}_{x/y/z}}) - n_{x/y/z}}{n_{x/y/z}} \right| \leq \left(1 + \frac{2}{|n_{x/y/z}|}\right) u^2 + O(u^3).
\end{equation}

\subsection{Accurate Summation of Squares}

 While summing squares is not prone to catastrophic cancellation, na\"ive summation still introduces errors at the level of machine epsilon, which are amplified by our subsequent calculations. To address this, we employ specialized accurate summation algorithms previously developed in the literature. In this section, we review these techniques and propose composable versions to support enhancing robustness of larger computations that use them as building blocks.

The \texttt{SumOfSquares} algorithm (Algorithm \ref{alg:sum_of_squares}) proposed in \cite{Graillat2015} accurately computes the squared norm of a vector by leveraging the \texttt{TwoProd} and \texttt{SumNonNeg} operations. The closely related \texttt{SumNonNeg} algorithm (Algorithm \ref{alg:sum_non_neg})  \cite{Graillat2015} sums non-negative numbers and achieves a relative error bound of \(3u^2\).

\begin{algorithm}[H]
\caption{\texttt{SumOfSquares}}
\label{alg:sum_of_squares}
\begin{algorithmic}[1]
\STATE \textbf{Input:} $\vec{x} \in \mathbb{F}^n$
\STATE Initialize $S^* \gets 0$, $s^* \gets 0$
\FOR{$j = 1$ to $n$}
    \STATE $[P, p] \gets \texttt{TwoProd}(x_j, x_j)$
    \STATE $[S^*, s^*] \gets \texttt{SumNonNeg}([S^*,s^*], [P,p])$
\ENDFOR
\STATE \textbf{Return:} $[S^*, s^*]$
\end{algorithmic}
\end{algorithm}

\begin{algorithm}[H]
\caption{\texttt{SumNonNeg}}
\label{alg:sum_non_neg}
\begin{algorithmic}[1]
\STATE \textbf{Input:} $[A, a], [B, b] \in \mathbb{F}^2$
\STATE $[H, h] \gets \texttt{TwoSum}(A, B)$
\STATE $c \gets \texttt{fl}(a + b)$
\STATE $d \gets \texttt{fl}(h + c)$
\STATE $[S^*, s^*] \gets \texttt{FastTwoSum}(H, d)$ 
\STATE \textbf{Return:} $[S^*, s^*]$
\end{algorithmic}
\end{algorithm}

\texttt{SumOfSquares} achieves high accuracy when summing exact inputs, but is not able to compensate for inaccuracy in the input vector.
In our intersection point computation, the sum-of-squares operation is used to compute squared norms of vectors obtained from our cross-product routine and must account for the estimated error components returned by that routine.
We therefore propose a composable version, \texttt{SumOfSquaresC}, accepting both an input vector $\bf x$ and a vector of componentwise compensation terms $\bf e$. We prove an error bound in \Cref{sec:sum_of_squares} of the Supplementary Material for the specific cases appearing in our application involving vectors in 2D and 3D with norm bounded by 1.

\begin{algorithm}[H]
\caption{Proposed \texttt{SumOfSquaresC}}
\label{alg:modified_sum_of_squares}
\begin{algorithmic}[1]
\STATE \textbf{Input:} Vectors $\vec{x}, \vec{e} \in \mathbb{F}^n$
\STATE $[S^*, s^*] \gets \texttt{SumOfSquares}(\vec{x})$
\STATE $R^* \gets \texttt{CompDot}(\vec{x}, \vec{e})$
\STATE $R \gets \texttt{FMA}(2, R^*, s^*)$
\STATE $[S, s] \gets \texttt{FastTwoSum}(S^*, R)$
\STATE \textbf{Return:} $[S, s]$
\end{algorithmic}
\end{algorithm}

The use of the fused multiply-add (\texttt{FMA}) operation is primarily for performance reasons and is not critical for achieving our accuracy goal. Additionally, \texttt{FastTwoSum} \cite{dekker1971floating} is employed instead of \texttt{TwoSum} to enhance efficiency. The final step of the algorithm, \([S, s] \gets \texttt{FastTwoSum}(S^*, R)\), is needed to ensure that the magnitude of the returned error compensation term is small ($|s| \le |S| u$) as assumed by our subsequent operations and analysis.


\section{Accurate intersection point calculation}
\label{sec:accurate_intersection_point_calculation}

Applying the EFT-based routines introduced in \Cref{sec:accurate_operators} to evaluate the formula in \Cref{eq:final_gca_constLat_intersection_point}, we obtain Algorithm \texttt{AccuX}, which efficiently computes intersection point coordinates \(\tilde{P}_x\) and \(\tilde{P}_y\) with high accuracy.

\begin{algorithm}[H]
\caption{Accurate Calculation of Intersection Points Using EFT (\texttt{AccuX})}
\label{alg:accurate_intersection_performance}
\begin{algorithmic}[1]
\STATE \textbf{Input:} Points $\mathbf{x}_1, \mathbf{x}_2$, and scalar $z_0$
\STATE $[\bar{n}_x, e_{\bar{n}_x}] \gets \texttt{AccuDOP}(\mathbf{x}_1[1], \mathbf{x}_2[2], \mathbf{x}_1[2], \mathbf{x}_2[1])$
\STATE $[\bar{n}_y, e_{\bar{n}_y}] \gets \texttt{AccuDOP}(\mathbf{x}_1[2], \mathbf{x}_2[0], \mathbf{x}_1[0], \mathbf{x}_2[2])$
\STATE $[\bar{n}_z, e_{\bar{n}_z}] \gets \texttt{AccuDOP}(\mathbf{x}_1[0], \mathbf{x}_2[1], \mathbf{x}_1[1], \mathbf{x}_2[0])$
\STATE $[S_2, s_2] \gets \texttt{SumOfSquaresC}(\{\bar{n}_x, \bar{n}_y\}, \{e_{\bar{n}_x}, e_{\bar{n}_y}\})$
\STATE $[S_3, s_3] \gets \texttt{SumOfSquaresC}(\{\bar{n}_x, \bar{n}_y, \bar{n}_z\}, \{e_{\bar{n}_x}, e_{\bar{n}_y}, e_{\bar{n}_z}\})$
\STATE $[C, e_C] \gets \texttt{TwoProd}(z_0, z_0)$
\STATE $[D, e_D] \gets \texttt{CompDotC}(\{S_3, S_3, s_3, s_3\}, \{C, e_C, C, e_C\})$
\STATE $[E,e_E] \gets \texttt{TwoSum}(S_2, -D)$ 
\STATE $e \gets \texttt{fl}(\texttt{fl}(s_2-e_D)+ e_E))$
\STATE $[\overline{s^2}, e_{\overline{s^2}}] \gets \texttt{FastTwoSum}(E, e)$ \quad \quad \COMMENT{Ensure $e_{\overline{s^2}}$ is small}

\STATE $[\bar{s}, e_{\bar{s}}] \gets \texttt{AccSqrt}(\overline{s^2}, e_{\overline{s^2}})$
\STATE $[F_{x}, e_{F_{x1}}] \gets \texttt{TwoProd}(\bar{n}_x, \bar{n}_z)$
\STATE $e_{F_{x2}} \gets \texttt{fl}(\bar{n}_x * e_{\bar{n}_z})$
\STATE $e_{F_{x3}} \gets \texttt{fl}(\bar{n}_z * e_{\bar{n}_x})$
\STATE $e_{F_{x}} \gets \texttt{fl}(\texttt{fl}(e_{F_{x1}} + e_{F_{x2}} )+ e_{F_{x3}})$
\STATE $x \gets \texttt{CompDot}(\{F_{x}, e_{F_{x}}, \bar{n}_y, e_{\bar{n}_y}, \bar{n}_y, e_{\bar{n}_y}\}, \{z_0, z_0,\bar{s}, \bar{s}, e_{\bar{s}}, e_{\bar{s}}\})$
\STATE $[F_{y}, e_{F_{y1}}] \gets \texttt{TwoProd}(\bar{n}_y, \bar{n}_z)$
\STATE $e_{F_{y2}} \gets \texttt{fl}(\bar{n}_y * e_{\bar{n}_z})$
\STATE $e_{F_{y3}} \gets \texttt{fl}\bar{n}_z * e_{\bar{n}_y})$
\STATE $e_{F_{y}} \gets \texttt{fl}(\texttt{fl}(e_{F_{y1}} + e_{F_{y2}}) + e_{F_{y3}})$
\STATE $y \gets \texttt{CompDot}(\{F_{y}, e_{F_{y}}, -\bar{n}_x, -e_{\bar{n}_x}, -\bar{n}_x, -e_{\bar{n}_x}\}, \{z_0,z_0, \bar{s}, \bar{s}, e_{\bar{s}}, e_{\bar{s}}\})$
\STATE $R_x \gets x / S_2$
\STATE $R_y \gets y / S_2$
\STATE \textbf{Return:} $\{-R_x, -R_y\}$
\end{algorithmic}
\end{algorithm}

We begin by using \texttt{AccuDOP} to obtain the normal vector components, each returned as a compensated pair. Next, we apply \texttt{SumOfSquaresC} to compute the terms $\norm{\n_{xy}}^2$ and $\|\mathbf{n}\|^2$, yielding the pairs $[S_2, s_2]$ and $[S_3, s_3]$, respectively. We calculate $z_0^2$ exactly as the pair $[C, e_C]$ and then accurately evaluate the product $\|\mathbf{n}\|^2 z_0^2$ as the pair $[D, e_D]$ using \texttt{CompDotC}.
We obtain our accurate evaluation of $s^2$ by running \texttt{TwoSum} on the leading-order parts, executing floating point addition of the error terms, and finally running \texttt{FastTwoSum} to ensure the relationship
$|e_{\overline{s^2}}| \le \overline{s^2} u$. This last operation is required  for the Taylor expansion used by \texttt{AccSqrt} \cite{Rump2023}  in the next step to produce an accurate approximation to $s$.
Next, \texttt{CompDot} is used to evaluate the expression $z_0 n_x n_z + s n_y$, which is computed as  $x \approx
(F_{x} + e_{F_{x}}) z_0 + (\bar{n}_y + e_{\bar{n}_y}) (\bar{s} + e_{\bar{s}})$. Finally, the computed value \( x \) is divided by \( S_2 \), the leading-order part of \( \norm{\n_{xy}}^2 \), to obtain the intersection coordinate $\tilde{P}_x$, and the same process is used to obtain $\tilde{P}_y$.

\subsection{Error Analysis of \texttt{AccuX}}
We first establish upper and lower bounds for coordinates $\tilde P_x$ and $\tilde P_y$ and then use these to obtain a bound on $\norm{\tilde {\bf P} - {\bf P}}$.
We do this by propagating the normal component error bounds (\Cref{eq:accu_normal_error}) through all subsequent computations.
In the following, we denote the lower and upper bounds an expression involving our computed quantities as \(\LB{\cdot}\) and \(\UB{\cdot}\), respectively.

The expressions in this derivation become quite involved, and
we rely on computer algebra (Wolfram Mathematica \cite{Mathematica}) to perform each step and simplify bounds to the extent possible. The remainder of this subsection describes the propagation rules we employ for each step of \Cref{alg:accurate_intersection_performance} while the next subsection presents and discusses our final error bound. The Mathematica notebook implementing these propagation rules has been provided in the supplemental material.

We begin by obtaining bounds for the computed versions of expressions $\norm{\n_{xy}}^2$, $\norm{\n}^2$, and $z_0^2$, which are \(S_2 + s_2\), \(S_3 + s_3\), and \(C + e_C\), respectively. These are derived and presented in the supplement (\Cref{sec:sum_of_squares}).
Next, we obtain bounds for \(D + e_D\), which corresponds to \(\|\mathbf{n}\|^2 z_0^2\):
\begin{align*}
\UB{D + e_D} &= \UB{S_3 + s_3}\UB{C + e_C}(1 + \UB{\varepsilon_{\texttt{CompDotC}}}) \\ 
 \LB{D + e_D }&= \LB{S_3 + s_3}\LB{C + e_C}
 - \UB{S_3 + s_3}\UB{C + e_C}\UB{\varepsilon_{\texttt{CompDotC}}}.
\end{align*}
Here, \(\UB{\varepsilon_{\texttt{CompDotC}}}\) refers to the error bound from \Cref{eq:dotfma_re_rel_err} evaluated using the upper bounds for all components of the vectors passed to \texttt{CompDotC} in this step.
Next, we obtain bounds for \(\overline{s^2} + e_{\overline{s^2}}\), the calculated approximation of $s^2$:
\begin{align*}
    \overline{s^2} + e_{\overline{s^2}} &= E+e \\
    &= E + (-e_D + s_2 + e_E) \left(1+\varepsilon_{\text{flSum}}\right) \\
    &= E + e_E - e_D + s_2 + (-e_D + s_2 + e_E) \varepsilon_{\text{flSum}} \\
    &= S_2 + s_2 - (D + e_D) + (-e_D + s_2 + e_E) \varepsilon_{\text{flSum}}, \\
    \UB{\overline{s^2} + e_{\overline{s^2}}} &= \UB{S_2 + s_2} - \LB{D + e_D} + \UB{(-e_D + s_2 + e_E) \varepsilon_{\text{flSum}}}, \\
    \LB{\overline{s^2} + e_{\overline{s^2}}} &= \LB{S_2 + s_2} - \UB{D + e_D} - \UB{(-e_D + s_2 + e_E) \varepsilon_{\text{flSum}}}.
\end{align*}
where \(\UB{(-e_D + s_2 + e_E) \varepsilon_{\text{flSum}}}\) is a bound for the compensation term given by:
$$
   \big(\UB{D + e_D}\,u + \UB{S_2 + s_2}\,u + \UB{|e_E|}\big)\UB{\varepsilon_{\text{flSum}}}, \\
$$
$$
\text{with} \quad 
\UB{|e_E|} = \big(\UB{S_2 + s_2}\left(1+u\right)-\LB{D + e_D}+\UB{D + e_D}\big)u.
$$
The term \(\varepsilon_{\text{flSum}}\) represents the relative error for summing an array of floating point numbers
with ordinary floating point operations.
According to Lemma 8.4 from \cite{Higham2002}, this relative error is bounded by
$\UB{\varepsilon_{\text{flSum}}} = \gamma_{n-1} \frac{\sum_{i=1}^{n} |a_i|}{\left| \sum_{i=1}^{n} a_i \right|}$
for the generic sum \(\sum_{i=1}^{n} a_i\),
regardless of the order of summation or potential underflow. 

These bounds are then propagated through the square root:
\begin{align*}
\UB{\tilde{s}} = \UB{\bar{s} + e_{\bar{s}}} &= \sqrt{\UB{\bar{s}^2 + e_{\bar{s}^2}}}(1 + \UB{\varepsilon_{\text{AccuSqrt}}}),
\\
\LB{\tilde{s}} = \LB{\bar{s} + e_{\bar{s}}} &= \sqrt{\LB{\bar{s}^2 + e_{\bar{s}^2}}} - \sqrt{\UB{\bar{s}^2 + e_{\bar{s}^2}}}\UB{\varepsilon_{\text{AccuSqrt}}},
\end{align*}
where \(\UB{\varepsilon_{\text{AccuSqrt}}} = \frac{25}{8} u^2\) is the relative error bound for the accurate square root operation derived in \cite{Rump2023}.

The calculated numerator subexpression $n_x n_z$, (\( F_{x} + e_{F_{x}} \)) can be expressed as:
\begin{align*}
    F_{x} + e_{F_{x}} &= F_{x} + \left(e_{F_{x_1}} + e_{F_{x_2}} + e_{F_{x_3}}\right)\left(1 + \varepsilon_{\text{flSum}}\right) \\
    &= F_{x} + e_{F_{x_1}} + e_{F_{x_2}} + e_{F_{x_3}} + \left(e_{F_{x_1}} + e_{F_{x_2}} + e_{F_{x_3}}\right) \varepsilon_{\text{flSum}} \\
    &= \bar{n}_x \bar{n}_z + \bar{n}_x e_{\bar{n}_z} \left(1+\delta_1\right) 
    + \bar{n}_z e_{\bar{n}_x} \left(1+\delta_2\right) 
    + \left(e_{F_{x_1}} + e_{F_{x_2}} + e_{F_{x_3}}\right) \varepsilon_{\text{flSum}} \\
    &= \tilde{n}_x \tilde{n}_z - e_{\bar{n}_x} e_{\bar{n}_z} 
    + \delta_1 \bar{n}_x e_{\bar{n}_z} 
    + \delta_2 \bar{n}_z e_{\bar{n}_x} 
    + \left(e_{F_{x_1}} + e_{F_{x_2}} + e_{F_{x_3}}\right) \varepsilon_{\text{flSum}},
\end{align*}

From this, we obtain the bounds:
\begin{align*}
    \UB{ F_{x} + e_{F_{x}}} &= \UB{\tilde{n}_x} \UB{\tilde{n}_z} + \UB{e_{\bar{n}_x}} \UB{e_{\bar{n}_z}} \\
    &\quad + u \left(\UB{\bar{n}_x} \UB{e_{\bar{n}_z}} 
    + \UB{\bar{n}_z} \UB{e_{\bar{n}_x}} \right) \\ 
    &\quad + \left(\UB{|e_{F_{x_1}}|} + \UB{|e_{F_{x_2}}|}  + \UB{|e_{F_{x_3}}|}\right) \UB{\varepsilon_{\text{flSum}}}, \\
    \LB{ F_{x} + e_{F_{x}}} &= \LB{\tilde{n}_x} \LB{\tilde{n}_z} - \UB{e_{\bar{n}_x}} \UB{e_{\bar{n}_z}} \\
    &\quad - u \left(\UB{\bar{n}_x} \UB{e_{\bar{n}_z}} 
    + \UB{\bar{n}_z} \UB{e_{\bar{n}_x}} \right) \\ 
    &\quad - \left(\UB{|e_{F_{x_1}}|} + \UB{|e_{F_{x_2}}|}  + \UB{|e_{F_{x_3}}|}\right) \UB{\varepsilon_{\text{flSum}}}.
\end{align*}

Finally, we obtain bounds for $R_x \defeq -\tilde{P}_x$, the intersection point coordinate computed before the final negation:
\begin{align*}
    \UB{R_x} &= \UB{\frac{x}{S_2}\left(1+\delta_3\right)} = \UB{ \frac{x}{S_2}}\left(1+\UB{\delta_3}\right)\\
     &=\frac{\left(\UB{F_{x}+e_{F_{x}}}z_0+\UB{\tilde{s}}\UB{\tilde{n}_y}\right)\left(1+u\right)}{\LB{S_2+s_2}-u\UB{S_2+s_2}}   \left(1+\UB{\delta_3}
    \right)+O(u^2),\\
    \LB{R_x} &= \LB{\frac{x}{S_2}\left(1+\delta_4\right)} = \LB{ \frac{x}{S_2}} - \UB{ \frac{x}{S_2}} \UB{\delta_4}\\
    &= \frac{\LB{F_{x}+e_{F_{x}}}z_0 + \LB{\tilde{s}}\LB{\tilde{n}_y} - \left(\UB{F_{x}+e_{F_{x}}}z_0+\UB{\tilde{s}}\UB{\tilde{n}_y}\right)u}{\UB{S_2+s_2}\left(1+u\right)}\\
    &\quad -  \UB{ \frac{x}{S_2}}\UB{\delta_4}+O(u^2),
\end{align*}

The bounds for $R_x$ were straightforward to derive because the corresponding expression
$\frac{z_0 n_x n_z + s n_y}{\norm{\n_{xy}}^2}$ involves adding together quantities known to be positive under our assumptions that $[n_x, n_y, n_z]$ are positive. However, deriving upper and lower bounds for $R_y = \frac{z_0 n_y n_z - s n_x}{\norm{\n_{xy}}^2}$ is more complicated: the computed numerator can be positive or negative as can its upper and lower bounds.
For example, when computing the upper bound of $R_y$, if the numerator's upper
bound is negative, it must be divided by the \emph{upper}, not lower, bound of
the denominator. There are eight possible combinations of signs to consider, although
two of them can be neglected (the upper bound must be positive when the
numerator itself is positive, and the lower bound must be negative when the
numerator is negative). We obtain these $R_y$ bounds separately for each possible case, and consider
each when obtaining the final error bound in the next subsection. use them to
compute an error bound for $\|\tilde{\vec{P}} - \vec{P}\|$, and pick the
largest.

\subsection{Error Analysis Results and Discussion}
Our error bound for $\|\tilde{\vec{P}} - \vec{P}\|$ is obtained by, for each component of $\vec{P}$, computing the maximum distance between its exact value and each of its upper and lower bounds
obtained in the previous subsection.
As in \Cref{sec:float_error}, we can divide our analysis into two regimes based
on the size of $s^2$. In regime (ii), where $|\flabserror{s^2}| \ll s^2$,
separation parameter $c < \infty$ can easily be chosen so that the
error term inflated by the square-root, $|\flabserror{s^2}|^{\frac{2c-1}{2c}}$, is higher order than the $O(u)$ error from the final division.
Specifically, since the \texttt{AccuX} evaluation ensures $|\flabserror{s^2}| = O(u^2)$, we simply need to pick $c$ such that:
$$
2\frac{2c-1}{2c} > 1 \quad \iff \quad c >1.
$$
This corresponds to ensuring $s > u^{1 + \tilde{c}}$ for some $\tilde c > 0$,
preventing the line of constant latitude from slicing excessively close to the apex
of the great circle arc.

Under this mild nondegeneracy assumption on the input,
we expect the leading-order term of $\flerror{P_x}$ and $\flerror{P_y}$ to be $3 u$,
accounting for three roundings to working precision (rounding the accurately-evaluated numerator and denominators of the final division, and then rounding the result).
And this is confirmed by our Mathematica analysis:
the regime (ii) assumption enables the use of a Taylor expansion to simplify
the square-root expressions otherwise appearing in the bound, leaving behind
\begin{equation}
\label{eq:final_rel_error_our1}
    \begin{aligned}
        \|\tilde{\vec{P}} - \vec{P}\| &\leq u \left\|\begin{bmatrix}
            \UB{\left|R_x - \frac{z_0 n_x n_z + s n_y}{\norm{\n_{xy}}^2}\right|} \\
            \UB{\left|R_y - \frac{z_0 n_y n_z - s n_x}{\norm{\n_{xy}}^2}\right|} \\
            0
        \end{bmatrix}\right\| \\
        &\leq u \left\|\begin{bmatrix}
            3P_x \\
            3P_y \\
            0
        \end{bmatrix}\right\| + \hot = 3 \sqrt{1 - z_0^2} u + \hot
    \end{aligned}
\end{equation}
The coefficient $3 \sqrt{1 - z_0^2}$ is no greater than $3$ for all inputs,
achieving our accuracy goal for this input regime.

In regime (i), we cannot use the Taylor expansion to simplify the Mathematica
bound, and our analysis in 
\Cref{eq:es_regime_i} indicates that the $O(u^2)$ error magnitude $|\flabserror{s^2}|$ gets amplified to $O(u)$ by the square root.
Ultimately this contributes atop the $3 u$ error from the final division an $O(u)$ term with
a complicated input-dependent coefficient reported in the supplementary Mathematica notebook.
This nevertheless represents a significant improvement over the accuracy of the direct floating-point implementation
(\Cref{eq:native_float_regime_i}).

We finally note that all bounds derived in this paper apply equally to the second intersection point:
$$
\small
\begin{aligned}
    P_x = \frac{z_0 n_x n_z - s n_y}{\norm{\n_{xy}}^2}, \quad
    P_y = \frac{z_0 n_y n_z + s n_x}{\norm{\n_{xy}}^2},
\end{aligned}
$$
as the formulas for the two points differ only by an interchange of the $x$ and $y$ axes.

\section{Experiments}
\label{sec:experiments}

To evaluate the performance and empirical accuracy of our proposed algorithm \texttt{AccuX}, we conduct several experiments comparing it to other implementations. Our benchmarks assess both computational speed and numerical accuracy, and we use Mathematica \cite{Mathematica} to generate ground-truth solutions and to compute the errors in our calculated results in arbitrary precision.

We specifically test the standard double-precision floating-point implementation of \cref{eq:final_gca_constLat_intersection_point},
as well as its evaluation using the GNU MPFR multiprecision library
\cite{MPFRWebsite, MPFR} at precision levels ranging from 16 to 32 digits, and
IEEE 754 quadruple-precision (binary128) from GNU \texttt{libquadmath}
\cite{libquadmath} via the \text{mp++} library \cite{mppp}. We also included a CGAL
implementation, which computes the exact intersection point in an algebraic
number type and then rounds it to floating-point using \texttt{CGAL::to\_double},
and the implementation of Krumm's method from Mirshak (\cite{mirshak2024intersections}).
We finally investigated the effects on performance and accuracy of using different mathematical formulas for the computations done in native floating point, quadruple precision, and MPFR.

The dataset for our accuracy experiments was generated as follows: Using GNU MPFR with \texttt{mpfr\_float\_1000} precision, we randomly sampled points on the unit sphere. The northern hemisphere was divided into bands based on latitude, with narrower bands near the equator and poles and wider bands at intermediate latitudes. Within each band, two random points were chosen as endpoints for the great circle arc, ensuring that their latitudes fell within the corresponding band and their longitudes spanned from \(0^\circ\) to \(360^\circ\). These latitude and longitude values were then converted into Cartesian coordinates. Finally, coordinate \(z_0\) of the line of constant lattitude was randomly sampled between the \(z\)-coordinates of the two endpoints of the great circle arc.

We also created a secondary dataset designed to include ill-conditioned scenarios close to the equator (arc endpoints within the lattitude band \((0^\circ, {10^{-4}}^\circ]\)). For this dataset, \(z_0\) was set close to the maximum latitude of the great circle arc: 
\(
z_0 = \sqrt{\frac{n_x^2 + n_y^2}{\|\mathbf{n}\|^2}} - r\), with a small positive random offset $r$ drawn from specified intervals dividing the range \(10^{-15}\) to \(10^{-2}\).

All implementations were coded in C++ and tested on an AMD Ryzen 9 5950X processor, featuring 16 cores and supporting the AVX2 SIMD set instruction (which features a 256-bit vector width). We used the same \texttt{gcc-13} compiler for all methods, chosen since GCC is needed for \texttt{libquadmath} support.

\subsection{Empirical Accuracy Analysis}



\begin{figure}[htbp]
    \centering
    \includegraphics[width=\textwidth]{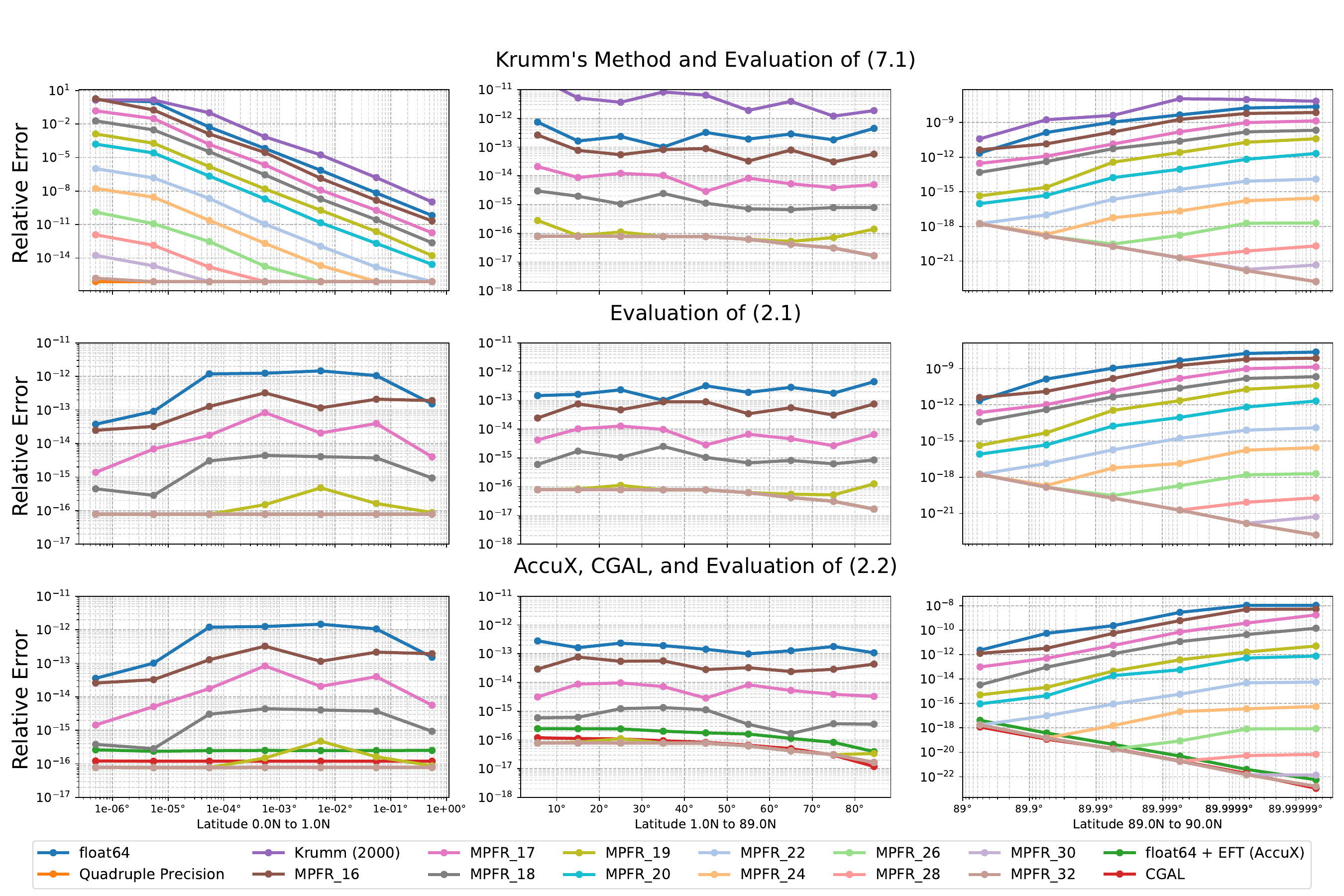}
\caption{Accuracy comparison of precision models (MPFR, quadruple-precision, and float64) for various intersection equations. The float64 + EFT (AccuX) represents our proposed \texttt{AccuX} algorithm. 
The first row shows accuracy using \cref{eq:gca_constlat_intersection_point_baseline} across different precision models, including results from Mirshak's Python code (2024) \cite{krumm2000,mirshak2024intersections}. The second row presents results with \cref{eq:gca_constlat_intersection_point_cdo}, and the third row uses our final proposed equation \cref{eq:final_gca_constLat_intersection_point}. Since AccuX is based on \cref{eq:final_gca_constLat_intersection_point}, it appears only in the third row. For comparison, CGAL's double-precision approximation is shown in the last row alongside AccuX.}
    \label{fig:comparison_methods}
\end{figure}

For each input in our primary dataset described above, we run all methods and compute relative errors with respect to the exact intersection points computed by Mathematica. To visually summarize the results, we group them by the latitude band from which the input arcs were generated and plot the worst (maximum) relative error encountered for each band. The resulting plots are shown in \Cref{fig:comparison_methods}, where the horizontal coordinate of each data point is obtained as the (arithmetic) midpoint of the latitude band.

The top row of \cref{fig:comparison_methods} includes the results from Mirshak's 2024 implementation of Krumm's method, which exhibits significant precision loss near the equator. To understand why this may happen, we introduce formula \cref{eq:gca_constlat_intersection_point_baseline} that exhibits a similar error trend and illustrates the sorts of numerical issues that can contribute to such precision loss around the equator:
\begin{equation}
\mathbf{P} = 
\begin{bmatrix}
    - \frac{1}{1 - \hat{n}_z^2}\left(z_0 \hat{n}_x \hat{n}_z \pm \hat{s}_* \hat{n}_y\right) \\
    - \frac{1}{1 - \hat{n}_z^2}\left(z_0 \hat{n}_y \hat{n}_z \mp \hat{s}_* \hat{n}_x\right) \\
    z_0
\end{bmatrix}, \quad \text{where } \hat{s}_* = \sqrt{1 - \hat{n}_z^2 - z_0^2}.
\label{eq:gca_constlat_intersection_point_baseline}
\end{equation}
We note that this new formula suffers from catastrophic cancellation when \(\hat{n}_z^2\) approaches 1, as the term \(1 - \hat{n}_z^2\) in the denominator approaches zero. Notably, even MPFR precision shows high relative errors with \Cref{eq:gca_constlat_intersection_point_baseline}. Using \Cref{eq:gca_constlat_intersection_point_cdo} (which avoids subtraction in the denominator), the relative error around the equator is greatly reduced.

The middle row of \cref{fig:comparison_methods} includes methods that evaluate \Cref{eq:gca_constlat_intersection_point_cdo} using different number types. This equation resolves the most severe cancellation issue of \Cref{eq:gca_constlat_intersection_point_baseline}, and the methods evaluating it achieve substantially more accurate results near the equator. This equation yields broadly similar accuracy behavior to \Cref{eq:final_gca_constLat_intersection_point}, since its extra normalization step is numerically stable. However, as we demonstrate \Cref{sec:performance}, these methods are substantially slower, and analyzing their error would be more challenging.

\begin{figure}[h!]
    \centering
    \includegraphics[width=\textwidth]{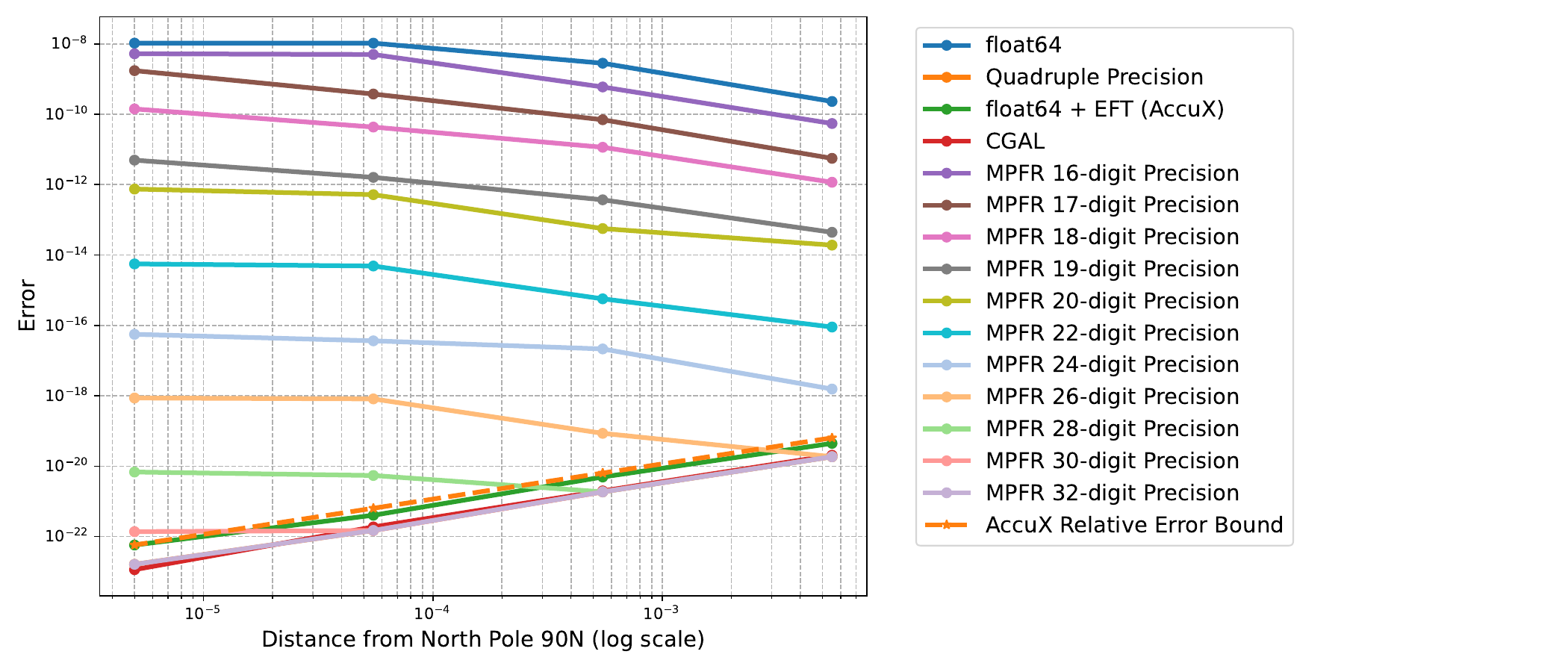}
\caption{Accuracy analysis of precision models near the North Pole (90N) using \cref{eq:final_gca_constLat_intersection_point} for MPFR, quadruple-precision, float64 (IEEE754 double precision), and CGAL’s double-precision approximation implementations. The float64 + EFT (AccuX) represents our proposed \texttt{AccuX} algorithm, demonstrating superior accuracy in the critical polar region compared to other methods.}

    \label{fig:pole_area_accuracy}
\end{figure}

Our method (\texttt{AccuX}) appears only in the final row of this figure, which also collects results from CGAL and methods based on evaluating \Cref{eq:final_gca_constLat_intersection_point}.
All methods perform reasonably well at the lower latitudes of this dataset (with native double-precision floating point only losing a few digits of precision), but the lower-precision methods exhibit significant error approaching the north pole. On the contrary, our method remains highly accurate, exhibiting similar error behavior to CGAL, quadruple precision, and MPFR32 (MPFR with 32 decimal digits). To better highlight the trends of each method as the intersection point approaches north pole, we provide \Cref{fig:pole_area_accuracy} which focuses on this part of our dataset.


\begin{figure}[H]
    \centering
    \includegraphics[width=\textwidth]{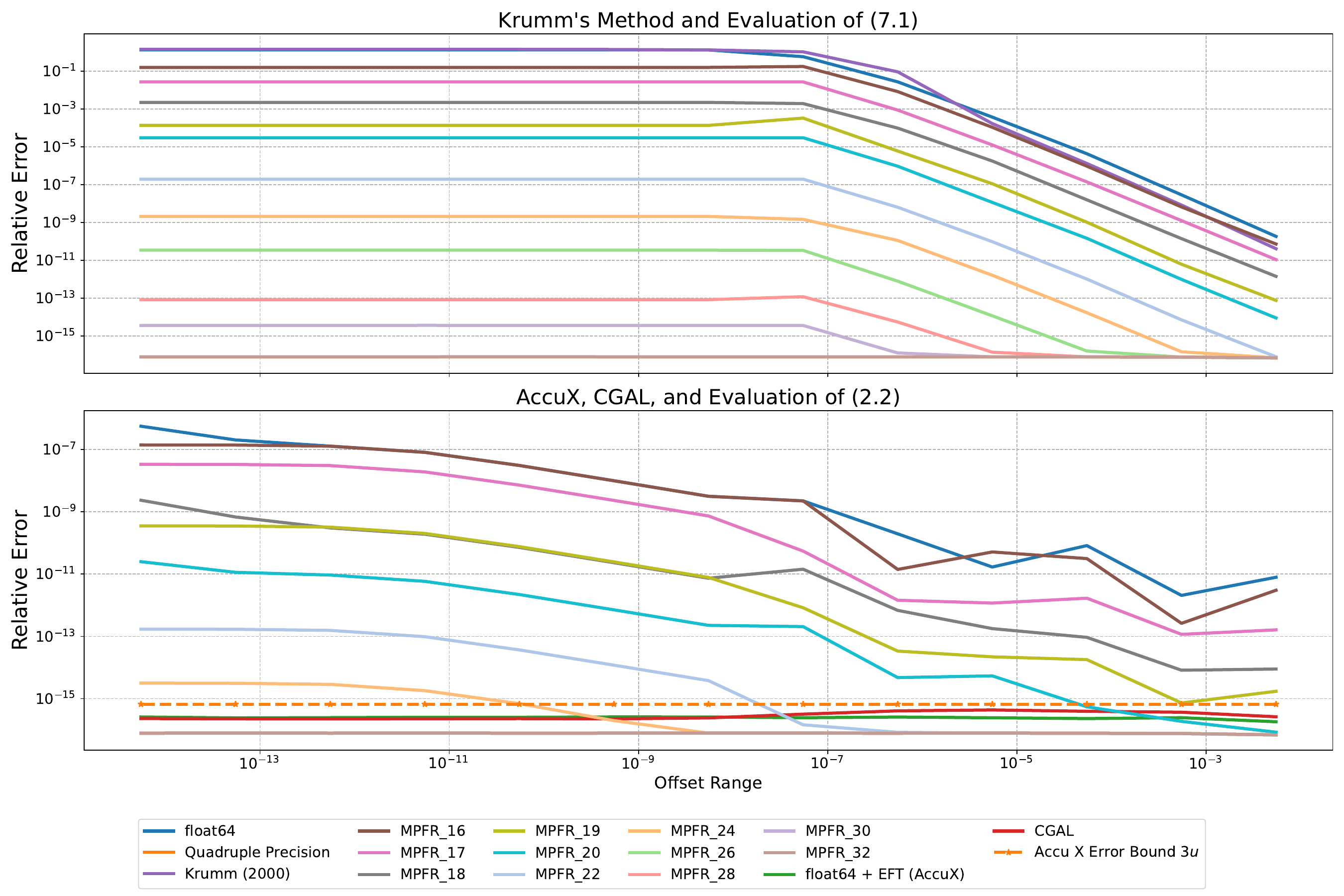}
\caption{Accuracy comparison of precision models (MPFR, quadruple-precision, and float64) for \cref{eq:gca_constlat_intersection_point_baseline} and \cref{eq:final_gca_constLat_intersection_point}, alongside Krumm’s (2000) method \cite{krumm2000}, in the calculation of ill-conditioned cases.}

    \label{fig:comparison_methods_extreme}
\end{figure}

Our primary dataset explored in \Cref{fig:comparison_methods} does not end up sampling a certain difficult class of inputs for which the intersection point calculation is ill-conditioned. This is the motivation for our secondary dataset containing shallow great circle arcs near the equator and lines of constant latitude slicing near the top of the arc.
We analyzed these inputs in the same way as our primary dataset, but we used the intervals from which offset $r$ was sampled, rather than the lattitude band, to group the inputs and define the horizontal axis.
These results are presented in \Cref{fig:comparison_methods_extreme}.
Only \texttt{AccuX}, CGAL, and high-precision MPFR maintain acceptable accuracy in this regime. Notably, Krumm's method exhibits relative errors exceeding 100\%.


To provide insight into how errors propagate through the evaluation of Equation \ref{eq:final_gca_constLat_intersection_point}, we computed the relative error of the intermediate quantities computed at each step of the calculation and plotted them in Figure \ref{fig:intermediate_results}. We did this both for \texttt{AccuX} and for a direct evaluation in double-precision and high-precision number types.

\subsection{Performance Analysis}
\label{sec:performance}
\begin{figure}[htbp]
    \centering
    \includegraphics[width=\linewidth]{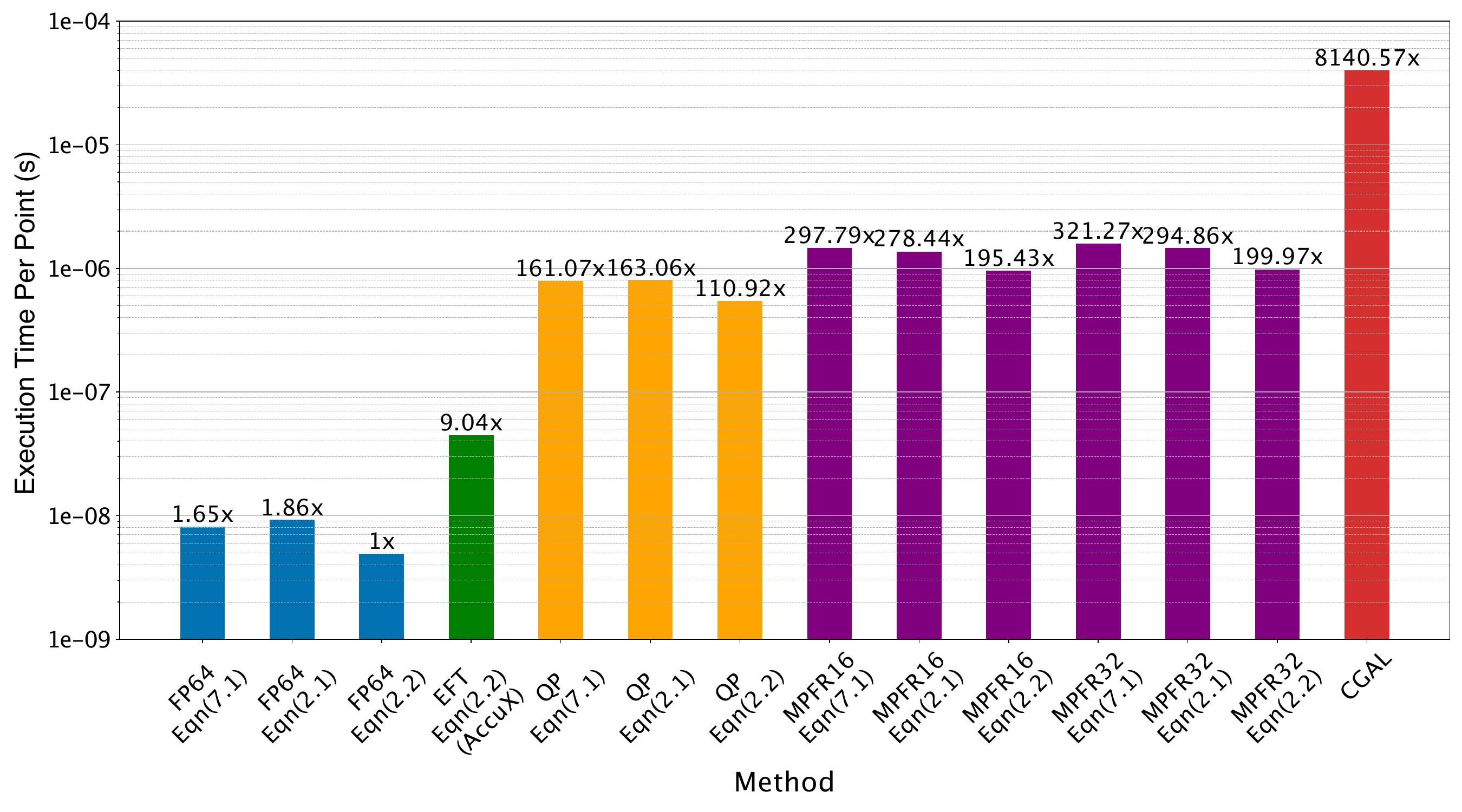} 
\caption{Performance comparison of different precision models, including MPFR, quadruple-precision (QP), and float64 (FP64) implementations, utilizing three mathematical formulas: \cref{eq:gca_constlat_intersection_point_baseline}, \cref{eq:gca_constlat_intersection_point_cdo}, and \cref{eq:final_gca_constLat_intersection_point}. FP64 + EFT (AccX) represents our proposed algorithm based on \cref{eq:final_gca_constLat_intersection_point}. CGAL presents the CGAL's double-precision approximation. Lastly, CGAL computes both intersection points, while all other methods compute only one, but with less than twice the workload due to shared intermediate results.}

    \label{fig:performance_results}
\end{figure}

\begin{figure}[htbp]
    \centering
    \includegraphics[width=\linewidth]{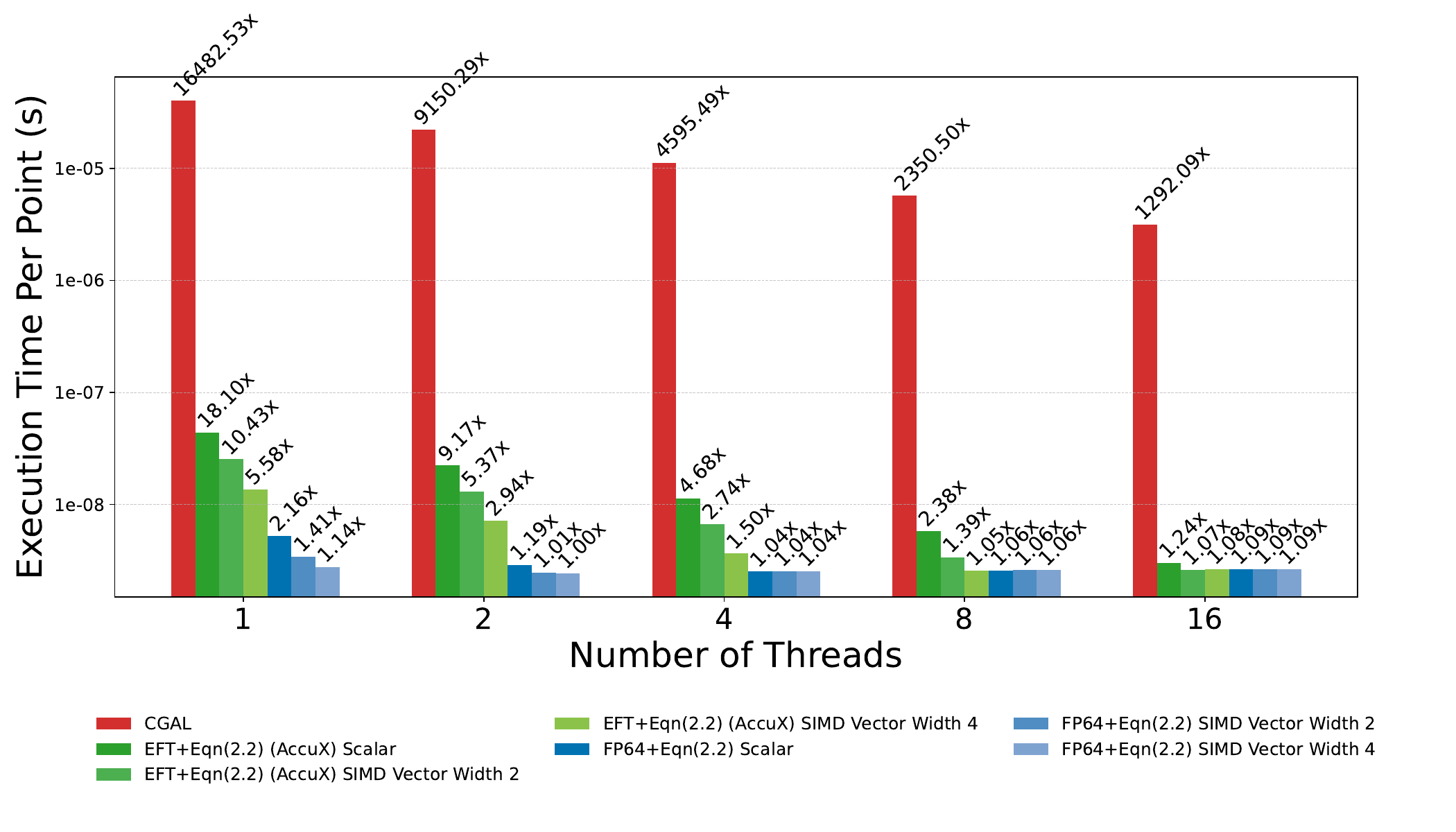}
\caption{Performance comparison of the \texttt{AccuX} algorithm (AccuX), an EFT method using \cref{eq:final_gca_constLat_intersection_point}, and the float64 (FP) implementation across various SIMD vector widths and thread counts. CGAL(double-precision approximation) results show different thread counts only. The y-axis represents execution time in seconds (log scale), and the x-axis shows the number of threads.}

    \label{fig:performance_comparison}
    \label{fig:performance_results_parallel}
\end{figure}

To analyze the trade-offs between numerical accuracy and computational cost, we benchmark each method on a dataset consisting of \(1000{,}000{,}000\) randomly generated pairs of intersecting arcs. Note that our intersection algorithm is branch-free, so its runtime is unaffected by the specific coordinates. Figure \ref{fig:performance_results} reports serial performance, while Figure \ref{fig:performance_results_parallel} shows scaling across different SIMD vector widths and thread counts. In both cases, \texttt{AccuX} achieves substantial speedups compared to all methods of comparable accuracy.

%
%
Our MPFR implementation is done using the Boost Multiprecision Library \cite{boost}, and we enable static stack allocation and expression templates for best performance.
 We timed both the 16-digit and 32-digit precision variants of this method, where the 16-digit precision approximates IEEE 754 double precision (\texttt{float64}), and 32-digit precision provides roughly quadruple precision.
For CGAL, we measured the time taken to compute the rounded intersection points (via \texttt{CGAL::to\_double}).
We note that CGAL computes all intersection point for a given input (up to two), while the other methods evaluate only the first solution listed in \Cref{eq:final_gca_constLat_intersection_point}.
However, considering the significant sharing of intermediate results between the two solutions, the overhead of computing all solutions would be negligible compared to increased time measured for CGAL.

Two key insights arise from the results:
First, the proposed formula \ref{eq:final_gca_constLat_intersection_point} demonstrates a substantial speedup over formula \ref{eq:gca_constlat_intersection_point_baseline} (by a factor of 1.65) and formula \ref{eq:gca_constlat_intersection_point_cdo} (by a factor of 1.86), largely due to avoiding the square-root operation needed for normalization. 
Second, leveraging the optimized formula, the \texttt{AccuX} algorithm consistently outperforms quadruple precision, both MPFR implementations, and CGAL, even in their highest-performance configurations. 
Our algorithm is approximately $9 \times$ slower than the pure floating-point implementation of \cref{eq:final_gca_constLat_intersection_point}, which is the fastest method overall. It is about $4.8\times$ slower than the floating-point implementations of both (\ref{eq:gca_constlat_intersection_point_baseline}) and (\ref{eq:gca_constlat_intersection_point_cdo}).
In contrast, when compared to quadruple precision (binary128), our method is roughly $12\times$ faster at evaluating (\ref{eq:final_gca_constLat_intersection_point}). Additionally, our algorithm is approximately $20\times$ faster than the MPFR library (whether 16 digit or 32) on the same formula. The CGAL implementation is the slowest method, being nearly four orders of magnitude slower than direct floating-point calculations.


For applications needing to compute the intersections of many input arcs, it makes sense to run the computations in parallel.
Since \texttt{AccuX} is built atop native floating point, it naturally supports SIMD (Single Instruction Multiple Data) parallelization using the CPU's vector instructions (for our system, this is \texttt{AVX2}).
This leads to a dramatic speedup, which is not possible with MPFR, quadruple-precision, or CGAL due to their
non-primitive number types. Furthermore, we can leverage multiple cores to parallelize the computation, and this is supported by all methods.

We implement vectorization using the Eigen library \cite{eigenweb}
and C++'s template functionality.
Specifically, \texttt{AccuX} and its underlying algorithms
are implemented in functions templated on the real number type to be used
for all inputs, outputs, and calculations.
This enables a standard scalar version of the code to be obtained by
instantiating the templates with type \texttt{double},
while a width-\texttt{N} SIMD version is obtained
using \texttt{Eigen::Array<double, N, 1>} as the number type.
We vectorize the baseline native floating point implementation
of \Cref{eq:final_gca_constLat_intersection_point} using the same
function template approach.

Our parallel performance experiments execute an outer loop across the $10^9 / \texttt{N}$ \emph{chunks} of
consecutive inputs in our benchmark dataset; this loop is parallelized
using \texttt{tbb::parallel\_for} from Intel's oneTBB
library \cite{intel_onetbb_2022_2_0}, which provides control
over the number of threads utilized.
In the vectorized case ($\texttt{N} > 1$), we load
the values of the \texttt{N} distinct inputs of the chunk into
array-valued variables holding $z_0$ and the individual coordinates of $\x_1$ and $\x_2$
before passing them to the appropriate vectorized
instantiation.
The code was compiled using the flags \texttt{-march=native -O3 -ffp-contract=off},
the last of which disables automatic generation of FMA instructions;
this ensures that FMAs are not introduced in places where they were not
explicitly intended by our algorithm to prevent
deviations of the rounding behavior from our numerical analysis.


Our benchmarking experiment results are presented in Figure \ref{fig:performance_results_parallel}, which compares the vectorized and parallelized implementations of the pure floating-point evaluation of (\ref{eq:final_gca_constLat_intersection_point}) with \texttt{AccuX} across various SIMD vector widths and thread counts. For reference, we also present CGAL's parallel performance with increasing thread counts.
The slowdown statistics annotated on each bar are computed relative to the fastest measured configuration
(native floating point, SIMD width $N=4$, and two threads).

Interestingly, 
the results show that our vectorized implementation of \texttt{AccuX} approaches and eventually closely matches the performance of the pure floating-point code
when both are run with increasing vector widths and thread counts.
This happens because the intersection point calculation is not computationally intensive, and the additional arithmetic required by \texttt{AccuX} ends up being hidden by memory I/O bottlenecks.
While the vectorized \texttt{AccuX} implementation continues to show parallel speedups
up to 8 threads before saturating memory bandwidth, trying to extract more performance from the direct floating-point implementation
by increasing beyond two threads fails (in fact leading to an slight slowdown, likely due to thread sheduling overhead).
We conclude that, for applications involving batch processing large sets of inputs on a parallel machine, the accurate results from \texttt{AccuX} can be obtained at
no compute-time overhead relative to a na\"ive floating point implementation.


As expected, the CGAL implementation remains orders-of-magnitude slower than \texttt{AccuX} after parallelization, even though it enjoys a near-perfect parallel speedup since it is not memory-bound.

\section{Conclusions}
\label{sec:conclusions}

We introduced a mathematical formula for calculating the intersection point between a geodesic arc and a constant latitude line on the unit sphere. This formula reduces floating-point rounding errors and enhances computational performance. Our error analysis of this formula revealed that implementing it in pure floating point arithmetic can lose half of the working precision for certain inputs. We address this by employing Error-Free Transformation (EFT) techniques to develop a substantially more accurate algorithm called \texttt{AccuX}.

Our algorithm guarantees that relative errors are bounded by \( 3\sqrt{1 - z_0^2}u \), ensuring a result with accuracy near machine-precision $u$. This accuracy rivals that of traditional high-precision methods such as MPFR, quadruple precision and CGAL, as confirmed empirically by our comprehensive accuracy benchmarks. At the same time, our method offers significantly better performance. In fact, although our algorithm performs more operations than a direct floating-point implementation, the compute time overhead vanishes once the codes are fully parallelized with vectorization and multithreading, eventually reaching the hardware's memory bandwidth limits. 


In summary, the \texttt{AccuX} algorithm delivers both high accuracy and computational efficiency, making it a scalable and robust solution for large-scale geodesic intersection problems. Though our work focuses on a specific geometric calculation,
its implications are broader: our algorithms and results demonstrate how EFT
techniques, when implemented with composability in mind, can be used to
solve non-trivial computations to high accuracy at much lower
overhead than with multiprecision number types.
They also show how, especially due the EFT approach's amenability to vectorization,
this overhead can vanish for large data-parallel applications.

\section{Future Work}
\label{sec:future_work}

While this work advances accuracy and efficiency in calculating intersection points between geodesic arcs and constant latitude lines, several areas for future exploration remain.
One direction is to apply the improved cross-product algorithm to other problems, such as the simpler one of computing intersections of two geodesic arcs. Another is to develop adaptive algorithms that use filters to dynamically switch between direct floating-point evaluation and our more accurate EFT-based calculations. However, introducing branching like this inhibits vectorization and, based on our performance experiments, would offer no speedups for large-scale batch processing.
Also, one could further investigate the relative error term at the $O(u^2)$ level in \Cref{eq:final_rel_error_our1} and develop a more detailed condition number for \texttt{AccuX}. This would provide deeper insight into the stability and precision limits of the algorithm, particularly in extreme cases where higher-order terms become significant.
Finally, because our algorithm is built on native floating-point operations, it is highly portable and could be used, for example, for GPU computation. More generally, it would be interesting to explore the use of EFT approaches to leverage the fast reduced-precision number types offered by modern GPUs while retaining sufficient accuracy.


\appendix

\section{Relative Errors of Intermediate Results} 
\label{sec:appendix2}
\Cref{sec:experiments} presented empirical accuracy experiments
evaluating each method based on the final result. However, analyzing
intermediate computational results can provide a deeper understanding of how
accuracy evolves throughout the computation and how errors propagate.
We specifically analyze the intermediate values for evaluating the intersection coordinate formula
$- \frac{1}{\norm{\n_{xy}}}^2 (z_0 n_x n_z + s n_y)$
on inputs from our primary dataset
using ordinary floating point, \texttt{AccuX}, and high-precision number types.

For \texttt{AccuX}, we recorded these results as pairs of doubles (\emph{i.e.},
the leading term and compensation term) and then computed the relative error
between their exact sum and the ground-truth quantity in Mathematica.
For example, to analyze the accuracy of $\norm{\tilde \n}^2$, we output the pair
$(S_3, s_3)$ computed by \Cref{alg:accurate_intersection_performance} and then
computed the relative error between $S_3 + s_3$ and $\norm{\n}^2$.
When using the high-precision number types from MPFR
and \texttt{libquadmath} via the \texttt{mp++ library}, we used these libraries' I/O routines to
export enough digits so that the values could be perfectly reconstructed in
Mathematica. 
Note that we did not include CGAL in this analysis since it computes
the exact intersection coordinate in an algebraic number type rather than by
evaluating the formula above.

\begin{figure}[htbp]
    \centering
    \includegraphics[width=\linewidth, trim=0 120 0 0, clip]{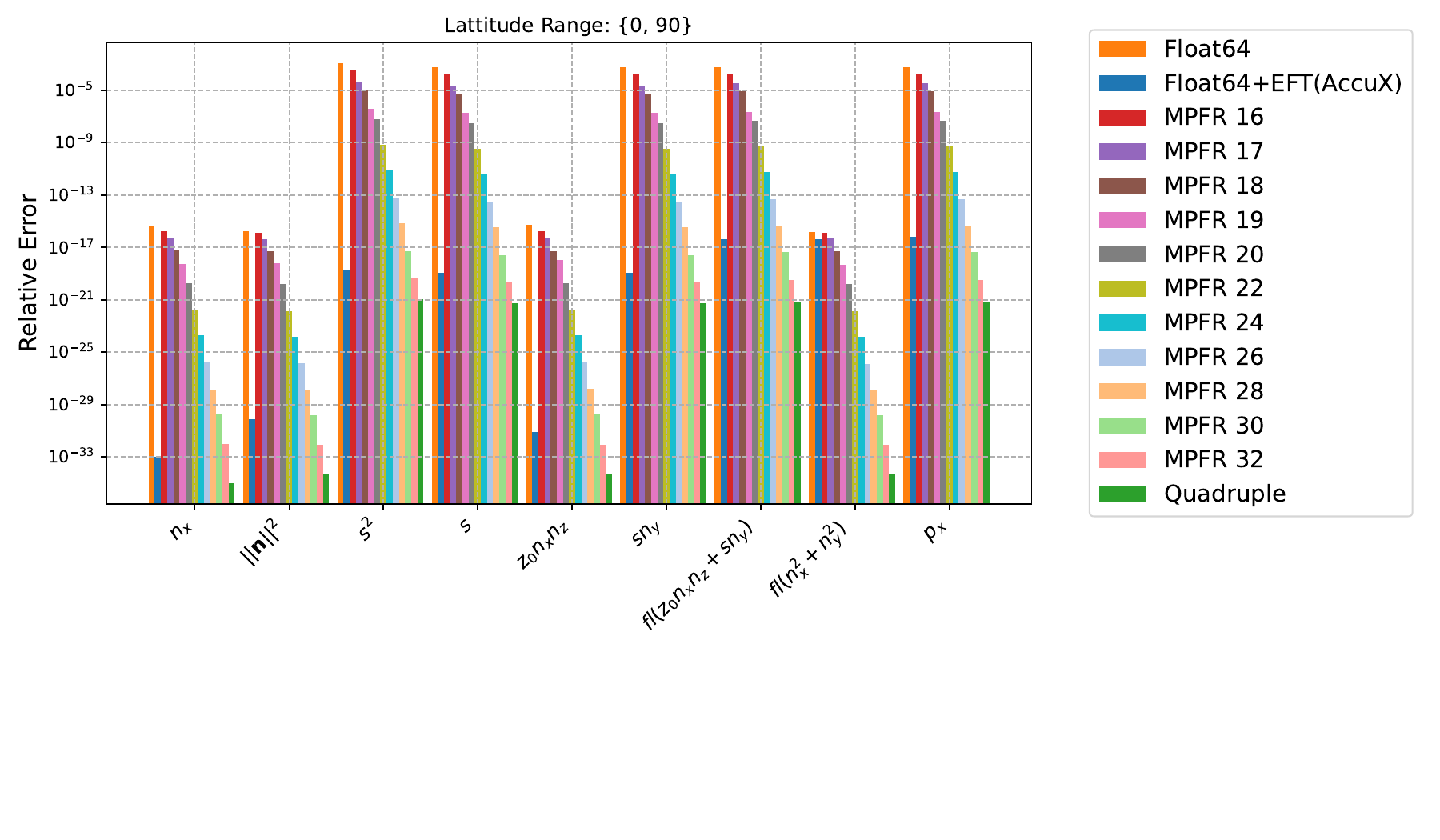}
    \caption{Relative error statistics for intermediate steps in evaluating the intersection coordinate formula.
             The reported errors have been averaged across all inputs in our primary dataset.
             The horizontal axis lists the intermediate computations that we analyzed, progressing from
             left to right, starting with an individual component of $\n$ and
             building up to the full coordinate (here labeled $p_x$).
            }
    \label{fig:intermediate_results}
\end{figure}

For each intermediate computation, we averaged its relative errors across the dataset to produce a single value plotted in \Cref{fig:intermediate_results}.
The resulting plot shows that our algorithm closely
matches the behavior of the higher-precision number types,
maintaining low relative errors throughout the computations.
Errors remain minimal for the cross product's components and squared norm $\norm{\n}^2$.
All intermediate results of \texttt{AccuX} exhibit relative errors below $u$ until the final
steps where the numerator and denominator are each rounded to
a single working-precision value (denoted as
$\text{fl}(z_0 n_x n_z + s n_y)$ and $\text{fl}(n_x^2 + n_y^2)$,
respectively) and then divided.
After these roundings, the relative error stabilizes close to $u$.

%% file: supplement_body.tex
\maketitle

This supplementary document first provides an error analysis of \Cref{eq:final_gca_constLat_intersection_point}
when implemented by direct translation into native floating point arithmetic.
Then it analyzes several of the building blocks of our accurate computation.

\section{Analysis of Direct Floating Point Implementation}
\label{sec:direct_fp_analysis}
In native floating point, each primitive arithmetic operation \( \star \in \{
    +,-,*,\divisionsymbol \} \) introduces rounding error that can be analyzed
using the following simple model \cite{Higham2002}:
\[
    x \flstar y = (x \star y) (1 + \delta) \quad \quad \abs{\delta} \le u,
\]
where \( \flstar \) denotes the floating point counterpart to \( \star \),
$\delta$ is the relative error introduced by rounding the exact result to the nearest floating point number,
and \( u \) is the machine epsilon.
While slightly tighter bounds can be obtained for each operation
\cite{JeannerodRump2017}, we use these classical bounds for
simplicity.

Our analysis in this section tracks how this roundoff error
propagates through all intermediate calculations of the intersection point $\P$
executed by a direct translation of the formulas into floating point.
To organize this analysis, we introduce the notation \( \flvalue{\subexpression}
\) to indicate the computed value of a given subexpression $\subexpression$ and
\( \flerror{\subexpression} \) to denote the relative error
\emph{measured with respect to the corresponding exact subexpression}.
Furthermore, we use $\flabserror{\subexpression}$ to denote the \emph{absolute} error of the subexpression.
These quantities are related by
$\flabserror{\subexpression} = \flvalue{\subexpression} - \subexpression$
and $\flabserror{\subexpression} = \subexpression\, \flerror{\subexpression}$.
This notation obeys the following recursive rules for combining subexpressions \(a \) and \( b \):
\begin{alignat*}{2}
    \flvalue{a \pm b} &= (\flvalue{a} \pm \flvalue{b}) (1 + \delta), \quad
    &&\flerror{a \pm b}
          = \frac{a (\flerror{a} + \delta) \pm b (\flerror{b} + \delta)}{a \pm b} + O(u^2),
\\
    \flvalue{a * b} &= (\flvalue{a} * \flvalue{b}) (1 + \delta), \quad
    &&\flerror{a*b} = \flerror{a} + \flerror{b} + \delta + O(u^2),
\\
    \flvalue{a \divisionsymbol b} &= \frac{\flvalue{a}}{\flvalue{b}} (1 + \delta), \quad
    &&\flerror{a \divisionsymbol b} = \flerror{a} - \flerror{b} + \delta + O(u^2).
\end{alignat*}
We use the rules above to obtain a bound for the absolute error $\norm{\flabserror{\P}}$, working through the evaluation tree of $\P$ from bottom to top.
At the leaves of this tree are the inputs \( \x_1 = (x_1, y_1, z_1) \), \( \x_2 = (x_2, y_2, z_2)\), and \( z_0 \), which we assume to be exact floating point numbers.
From these inputs, we compute the normal na\"ively:
\begin{align*}
    \flabserror{n_x} &= \big((y_1 z_2 (1 + \delta_1) - z_1 y_2 (1 + \delta_2)\big) (1 + \delta_3) - n_x
    \\
             &= \big((y_1 z_2 (\delta_1 + \delta_3) - z_1 y_2 (\delta_2 + \delta_3)\big) + O(u^2),
\\
    |\flabserror{n_x}| &\le (|y_1 z_2| + |z_1 y_2|)2 u + O(u^2) \le 2 u + O(u^2),
\end{align*}
where we used the bound $(|y_1 z_2| + |z_1 y_2|) \le 1$ justified in the main paper.
Note that the subscripted symbols $\delta_i$ refer to the relative error introduced by specific operations,
but ultimately we use the common bound $|\delta_i| \le u$ to remove subscripts.
We conclude:
\begin{align*}
    |\flabserror{n_x}| &\le 2 u + O(u^2), \quad |e_{n_y}| \le 2 u + O(u^2), \quad |e_{n_z}| \le 2 u + O(u^2), 
\\
    |\varepsilon_{n_x}| &\le \frac{2 u}{|n_x|} + O(u^2), \quad |\varepsilon_{n_y}| \le \frac{2 u}{|n_y|} + O(u^2), \quad |\varepsilon_{n_z}| \le \frac{2 u}{|n_z|} + O(u^2).
\end{align*}
Next, we introduce vector $\n_{xy} = [n_x, n_y]^\top$ and consider quantities $\norm{\n_{xy}}^2$ and $\norm{\n}^2$:
\begin{align*}
    \flerror{n_x^2} &= 2 \flerror{n_x} + \delta_4 + O(u^2), \; \flerror{n_y^2} = 2 \flerror{n_y} + \delta_5 + O(u^2), \; \flerror{n_z^2} = 2 \flerror{n_z} + \delta_6 + O(u^2),
    \\
    |\flerror{n_x^2}| &\le \frac{4 u}{|n_x|} + u + O(u^2), \; |\flerror{n_y^2}| \le \frac{4 u}{|n_y|} + u + O(u^2), \; |\flerror{n_z^2}| \le \frac{4 u}{|n_z|} + u + O(u^2),
    \\
    \flerror{\norm{\n_{xy}}^2} &= \frac{n_x^2 \flerror{n_x^2} + n_y^2 \flerror{n_y^2}}{\norm{\n_{xy}}^2} + \delta_7 + O(u^2),
    \\
    |\flerror{\norm{\n_{xy}}^2}| &\le \frac{4 u (|n_x| + |n_y|)}{\norm{\n_{xy}}^2} + 2 u + O(u^2)
                = \frac{4 \norm{\n_{xy}}_1}{\norm{\n_{xy}}^2} u + 2 u + O(u^2)
    \\
    \flerror{\norm{\n}^2} &= \frac{n_x^2 (\flerror{n_x^2} + \delta_7) + n_y^2 (\flerror{n_y^2} + \delta_7) + n_z^2 \flerror{n_z^2}}{\norm{\n}^2} + \delta_8 + O(u^2),
    \\
    \left|\flerror{\norm{\n}^2}\right| &\le \frac{4 u (|n_x| + |n_y| + |n_z|)}{\norm{\n}^2} + 3 u + O(u^2)
                = \frac{4 \norm{\n}_1}{\norm{\n}^2} u + 3 u + O(u^2). 
\end{align*}
Using these, we analyze the accuracy of $s^2 = \norm{\n_{xy}}^2 - \norm{\n}^2 z_0^2$:
\begin{align*}
    \left|\flerror{\norm{\n}^2 z_0^2}\right| &\le \left|\flerror{\norm{\n}^2}\right| + 2 u + O(u^2) \le \frac{4 \norm{\n}_1}{\norm{\n}^2} u + 5 u + O(u^2)
    \\
    \flabserror{s^2} &= \norm{\n_{xy}}^2 \flerror{\norm{\n_{xy}}^2} - \norm{\n}^2 z_0^2 \flerror{\norm{\n}^2 z_0^2} + s^2 \delta_9 + O(u^2)
    \\
    \left|\flabserror{s^2}\right| &\le  4 (\norm{\n_{xy}}_1 + z_0^2 \norm{\n}_1) u + 2 \norm{\n_{xy}}^2 u + 5 z_0^2 \norm{\n}^2 u + s^2 u + O(u^2)
    \\
                                  &\le  \Big(4 (\norm{\n_{xy}}_1 + z_0^2 \norm{\n}_1) + 7 \norm{\n_{xy}}^2 + s^2\Big) u + O(u^2),
\end{align*}
where we used $s^2 \ge 0 \implies z_0^2 \norm{\n}^2 \le \norm{\n_{xy}}^2$ to merge terms $z_0^2 \norm{\n}^2$ and $\norm{\n_{xy}}^2$.
We apply this inequality again along with $|z_0| \norm{\n}_1 \le \sqrt{3} |z_0| \norm{\n} \le \sqrt{3} \norm{\n_{xy}}$ and $\norm{\n_{xy}}_1 \le \sqrt{2} \norm{\n_{xy}}$ to further collect terms:
\begin{align*}
        \flabserror{s^2}
                                  &\le  \Big(4 (\sqrt{2} \norm{\n_{xy}} + \sqrt{3} |z_0| \norm{\n_{xy}}) + 7 \norm{\n_{xy}}^2 + s^2\Big) u + O(u^2)
    \\
                                  &\le 21 \norm{\n_{xy}} u + O(u^2).
\end{align*}
This loosens the bound by at worst a constant factor (at most
${21 / 4 = 5.25}$). It does not introduce any new expressions with potentially
worse growth behavior with respect to the input.

Next, we consider the square root used to obtain $s$ from $s^2$.
The behavior of the absolute error in $s$,
$$
    \flabserror{s} = s - \sqrt{s^2 + \flabserror{s^2}}(1 + \delta_{\sqrt{ }}),
$$
differs depending on whether $|\flabserror{s^2}| \ll s^2$.
We thus consider two regimes: (i) where $|\flabserror{s^2}| \approx s^2$ (but still $|\flabserror{s^2}| \le s^2$); and
(ii) where $|\flabserror{s^2}| \ll s^2$.
We therefore assume $\frac{|\flabserror{s^2}|}{s^2} \le 1$ and obtain the bound:
\begin{align}
    \flabserror{s} &= s\left(1 - \sqrt{1 + \frac{\flabserror{s^2}}{s^2}} + \sqrt{1 + \frac{\flabserror{s^2}}{s^2}} \delta_{\sqrt{ }}\right),
    \notag
    \\
    |\flabserror{s}| &\le s \left|1 - \sqrt{1 + \frac{\flabserror{s^2}}{s^2}}\right| + \sqrt{2} s |\delta_{\sqrt{ }}|
    \notag
    \\
                     &\le \frac{\left|\flabserror{s^2}\right|}{s} + \sqrt{2} s |\delta_{\sqrt{ }}|,
    \label{eq:es_bound}
\end{align}
where in the last step we used the upper bound of the inequality:
$$
(\sqrt{2} - 1) |x| \le |1 - \sqrt{1 + x}| \le |x| \quad \text{for} \quad x \in [-1, 1].
$$
The lower bound of this inequality indicates this last step loosens the bound by at most a factor of $1 + \sqrt{2}$. Furthermore,
when $x \to -1$ this upper bound is tight, and when $|x| \to 0$ it exceeds first-order Taylor approximation $|x| / 2$ by only a factor of $2$.

To analyze regime (i), which happens when there is significant cancellation error in evaluating (a tiny) $s^2$, we consider the limit of this bound as $\frac{|\flabserror{s^2}|}{s^2} \to 1$ from below:
\begin{equation}
\label{eq:es_regime_i}
    \UB{\left|\flabserror{s}\right|} = \sqrt{\left|\flabserror{s^2}\right|} \left( \frac{\sqrt{\left|\flabserror{s^2}\right|}}{s} + \sqrt{2} \frac{s}{\sqrt{\left|\flabserror{s^2}\right|}} |\delta_{\sqrt{ }}| \right)
            \to \sqrt{\left|\flabserror{s^2}\right|} + \hot
\end{equation}
The same limit process for the lower bound reveals the regime (i) error is exactly the same order as $\sqrt{\left|\flabserror{s^2}\right|}$:
$$
(\sqrt{2} - 1) \sqrt{\left|\flabserror{s^2}\right|} + \hot \le \left|\flabserror{s}\right| \le \sqrt{\left|\flabserror{s^2}\right|} + \hot
$$
In other words \emph{effectively half of the working precision is lost} for this class of inputs.

We note that when $|\flabserror{s^2}| > s^2$, which is not as relevant since
the intersection is likely to be missed ($\flvalue{s^2} < 0$), we can easily obtain the bound:
$$
\left|\flabserror{s}\right| \le |s| + \sqrt{{\flabserror{s^2}}}\sqrt{\frac{s^2}{|\flabserror{s^2}|} + 1} + \hot \le (1 + \sqrt{2}) \sqrt{|\flabserror{s^2}|} + \hot,
$$
also indicating a degradation to half the working precision.

For regime (ii),
we formalize our meaning of $|\flabserror{s^2}| \ll s^2$ by assuming
$|\flabserror{s^2}|^{1/c} < s^2$ for some $c > 1$ (\ie, a given level of separation between $|\flabserror{s^2}|$ and $s^2$).
This enables removing $s$ from the denominator of \Cref{eq:es_bound} in exchange for
lowering the order of the error expression $|\flabserror{s^2}|$ by factor $\frac{2c - 1}{2c}$:
\begin{equation}
    \label{eq:es_regime_ii_c}
    \left|\flabserror{s}\right| \le \frac{|\flabserror{s^2}|}{{|\flabserror{s^2}|}^{1/(2c)}} + s \sqrt{2} |\delta_{\sqrt{}}|
                             = {|\flabserror{s^2}|}^{\frac{2c - 1}{2c}} + s \sqrt{2} |\delta_{\sqrt{}}|
\end{equation}
As
$c \to 1$ (no separation), the order is lowered by a factor of $\frac{1}{2}$,
matching our findings for regime (i).
As $c \to \infty$ (infinite separation), no order lowering occurs.
Provided $c < \infty$ and $|\delta_{\sqrt{}}| = O(|\flabserror{s^2}|)$, the second term in \Cref{eq:es_regime_ii_c} is higher order
and can be neglected.

Next, we consider the numerator calculation:
\begin{align*}
    \flabserror{s n_y} &= s n_y \flerror{s n_y} = s n_y (\flerror{s} + \flerror{n_y} + \delta_{11} + \hot) = n_y e_s + s \flabserror{n_y} + s n_y \delta_{11} + \hot,
    \\
    \left|\flabserror{s n_y}\right| &\le |n_y| |e_s| + s(2 u + |n_y| u) + \hot,
    \\
    \flabserror{n_x n_z z_0} &= n_x n_z z_0 (\flerror{n_x} + \flerror{n_z} + \delta_{12} + \delta_{13} + O(u^2)),
    \\
    |\flabserror{n_x n_z z_0}| &\le 2 u |z_0| (|n_z| + |n_x| + |n_x n_z|) + O(u^2) \le 4 \norm{\n_{xy}} + O(u^2),
\end{align*}
where the last inequality used $|z_0|(|n_z| + |n_x| + |n_x n_z|) \le (\sqrt{2} + 1/2) |z_0| \norm{\n} < 2 \norm{\n_{xy}}$.
\begin{align*}
    \flabserror{(n_x n_z z_0 - s n_y)} &= \flabserror{n_x n_z z_0} - \flabserror{s n_y} + (n_x n_z z_0 - s n_y) \delta_{14} + \hot,
    \\
    |\flabserror{(n_x n_z z_0 - s n_y)}| &\le |\flabserror{n_x n_z z_0}| + |\flabserror{s n_y}| + |n_x n_z z_0 - s n_y| u + \hot
\end{align*}
Finally, we tackle the full coordinate evaluation:
\begin{align*}
    \flabserror{P_x} &= P_x \flerror{P_x} = P_x \left(\flerror{(n_x n_z z_0 - s n_y)} - \flerror{\norm{\n_{xy}}^2} + \delta_{15} + \hot\right)
    \\
                                         & = \frac{\flabserror{(n_x n_z z_0 - s n_y)}}{\norm{\n_{xy}}^2} - P_x \flerror{\norm{\n_{xy}}^2} + P_x \delta_{15} + \hot,
\\
    |\flabserror{P_x}| &\le \frac{\left|\flabserror{(n_x n_z z_0 - s n_y)}\right|}{\norm{\n_{xy}}^2} + |P_x| \flerror{\norm{\n_{xy}}^2} + |P_x| u + \hot
\\
                       &\le \frac{|\flabserror{n_x n_z z_0}| + |\flabserror{s n_y}|}{\norm{\n_{xy}}^2} + |P_x| \flerror{\norm{\n_{xy}}^2} + |P_x| 2u + \hot
\\
                       &\le \frac{|\flabserror{n_x n_z z_0}| + |\flabserror{s n_y}|}{\norm{\n_{xy}}^2} + |P_x| \left(\frac{4 \norm{\n_{xy}}_1}{\norm{\n_{xy}}^2} u + 4 u \right) + \hot
\\
                       &\le \frac{4 u \norm{\n_{xy}} + |n_y| |e_s| + s(2 u + |n_y| u)}{\norm{\n_{xy}}^2} + |P_x| \left(\frac{4 \sqrt{2} \norm{\n_{xy}}}{\norm{\n_{xy}}^2} u + 4 u \right) + \hot
\end{align*}
To obtain a simple, albeit looser, bound we collect terms using $s \le \norm{\n_{xy}}$:
\begin{align*}
    |\flabserror{P_x}| &\le \frac{|n_y| |e_s|}{\norm{\n_{xy}}^2} + \left. \frac{1}{\norm{\n_{xy}}} \middle(4 + 4 \sqrt{2} |P_x| + 2 + |n_y| \right) u + 4 |P_x| u + \hot
    \\
                       &\le \frac{|n_y| |e_s|}{\norm{\n_{xy}}^2} + \frac{7 + 4 \sqrt{2}}{\norm{\n_{xy}}} u + 4 |P_x| u + \hot
\end{align*}

\subsection{Regime (i) Point Error Bound}
In regime (i), the $x$ coordinate error is bounded by:
\begin{align*}
    |\flabserror{P_x}| \le \frac{|n_y| \sqrt{21 \norm{\n_{xy}}}}{\norm{\n_{xy}}^2} \sqrt{u} + \hot
\end{align*}
and likewise:
\begin{align*}
    |\flabserror{P_y}| \le \frac{|n_x| \sqrt{21 \norm{\n_{xy}}}}{\norm{\n_{xy}}^2} \sqrt{u} + \hot
\end{align*}
This makes the error bound for the full point:
\begin{align*}
    \norm{\flabserror{\P}} \le
            \sqrt{|\flabserror{P_x}|^2 + |\flabserror{P_y}|^2}
             = \norm{\n_{xy}} \frac{\sqrt{21 \norm{\n_{xy}}}}{\norm{\n_{xy}}^2} \sqrt{u} + \hot = \sqrt{\frac{21}{\norm{\n_{xy}}} u} + \hot
\end{align*}

\subsection{Regime (ii) Point Error Bound}
In regime (ii), when $1 < c < \infty$ we obtain the slightly more complicated bounds:
\begin{align*}
    |\flabserror{P_x}| &\le \frac{|n_y| (21 \norm{\n_{xy}} u)^\frac{2c - 1}{2c}}{\norm{\n_{xy}}^2} + \hot
    \\
    |\flabserror{P_y}| &\le \frac{|n_x| (21 \norm{\n_{xy}} u)^\frac{2c - 1}{2c}}{\norm{\n_{xy}}^2} + \hot
    \\
    \norm{\flabserror{\P}} &\le
            \sqrt{|\flabserror{P_x}|^2 + |\flabserror{P_y}|^2}
             = \frac{(21 \norm{\n_{xy}} u)^\frac{2c - 1}{2c}}{\norm{\n_{xy}}} + \hot
             \\
             &\le \frac{21}{\sqrt{\norm{\n_{xy}}}^{\frac{1}{c}}} u^\frac{2c - 1}{2c} + \hot
\end{align*}

\subsection{Error Bound in Terms of $s$}
We can also use \Cref{eq:es_bound} to obtain a bound on $\norm{\flabserror{\P}}$
in terms of $s$, without any assumptions of separation. Beginning with the components, and assuming $|\delta_{\sqrt{}}| \le u$
\begin{align*}
    |\flabserror{P_x}| &\le \frac{|n_y| \left(\left|\flabserror{s^2}\right| + s^2 \sqrt{2} |\delta_{\sqrt{}}|\right)}{s\norm{\n_{xy}}^2} + \frac{7 + 4 \sqrt{2}}{\norm{\n_{xy}}} u + 4 |P_x| u + \hot
    \\
                       &\le \frac{|n_y| \left(21 \norm{\n_{xy}} u\right)}{s\norm{\n_{xy}}^2} + \frac{11 + 5 \sqrt{2}}{\norm{\n_{xy}}} u + \hot
    \\
                       &\le \frac{21 |n_y|}{s\norm{\n_{xy}}}u + \frac{19}{\norm{\n_{xy}}} u + \hot
    \\
    |\flabserror{P_y}| &\le \frac{21 |n_x|}{s\norm{\n_{xy}}}u + \frac{19}{\norm{\n_{xy}}} u + \hot
\end{align*}
We then use the fact that:
\begin{align*}
\norm{\matrix{\frac{a |n_x|}{s \norm{\n_{xy}}} u + \frac{b}{\norm{\n_{xy}}} u + \hot
           \\ \frac{a |n_y|}{s \norm{\n_{xy}}} u + \frac{b}{\norm{\n_{xy}}} u + \hot
     }}^2
    &=
\left(
    \frac{a^2 \norm{\n_{xy}}^2}{s^2 \norm{\n_{xy}}^2}
    + \frac{2 a b \norm{\n_{xy}}_1}{s \norm{\n_{xy}}^2}
    + \frac{2 b^2}{\norm{\n_{xy}}^2} 
\right) u^2 + \hot
\\
    &\le 
\left(
    \frac{a^2}{s^2}
    + \frac{2 a b \sqrt{2}\norm{\n_{xy}}}{s \norm{\n_{xy}}^2}
    + \left(\frac{\sqrt{2} b}{\norm{\n_{xy}}} \right)^2
\right) u^2 + \hot
\\
    &=
\left(
    \frac{a}{s}
    +
    \frac{\sqrt{2} b}{\norm{\n_{xy}}}
\right)^2 u^2 + \hot
\end{align*}
to conclude:
$$
\norm{\flabserror{\P}} \le \left(\frac{21}{s} + \frac{19 \sqrt{2}}{\norm{\n_{xy}}}\right) u + \hot
$$
This bound indicates that error can blow up if $s \to 0$ (the line of constant latitude slices close to the apex of the great circle arc)
or if $\norm{\n_{xy}} \to 0$ (the great circle arc either approaches the equator or has its defining endpoints $\x_1$ and $\x_2$ approach each other).
Since $s \le \norm{\n_{xy}}$, we can summarize this bound as:
\begin{equation}
   \norm{\flabserror{\P}} \le \frac{48 u}{s} + \hot 
\end{equation}

\subsection{Impact of Using Kahan's $2 \times 2$ Determinant Algorithm}
We note that using Kahan's Algorithm to evaluate $\n$
will substantially improve accuracy by removing 
factors of $\norm{\n_{xy}}^{-1}$ from the bounds,
but it results in the same poor asymptotic behavior for $\norm{\flabserror{\P}}$,
yielding a leading-order term proportional to $\sqrt{u}$
in regime (i) and $u/s$ in regime (ii).

\section{Improved Accurate Sum of Squares Error Analysis}
\label{sec:sum_of_squares}
We have presented our proposed \texttt{SumOfSquaresC} algorithm in \Cref{sec:accurate_operators} and now provide a bound for its error.
In Graillat, Stef et al. (2015), the relative error bounds for the closely related \texttt{SumNonNeg} and \texttt{SumOfSquares} algorithms are derived. Specifically, for a floating point array $[x_1, x_2,..., x_n]$ with $\sigma = \sum_{i=1}^nx_i^2$ they proved that the result $[S^*,s^*]$ returned by \texttt{SumOfSquares} obeys 
$$
\frac{|(S^*+s^*)-\sigma|}{\sigma} \leq \frac{u}{8},
\quad
\text{when}
\quad
n \;<\; \bigl((24 + u)\,u \bigr)^{-1}.
$$
In addition \texttt{SumNonNeg} has an established relative error bound $\varepsilon_{\text{SNN}}\leq 3u^2$, provided the exact result is between $\frac{F_{\text{small}}}{u^2} $ and $ \frac{F_{\text{large}}}{2}$ where $F_{\text{small}}$ and $F_{\text{large}}$ are the smallest and largest positive normalized floating-point numbers, respectively \cite{Graillat2015}.

The bound above for \texttt{SumOfSquares} is intended for applications computing the squared norm of long vectors (large $n$), whereas the vectors in our application, $\n_{xy}$ and $\n$, are only of length 2 or 3. In these small cases, we can obtain a tighter error bound for it as follows. We denote by $\bigl[P_j, p_j\bigr]$ the pair computed in iteration $j$ of \Cref{alg:sum_of_squares} when computing the squared norm of $\bar{\n} = [\bar{n}_x, \bar{n}_y, \bar{n}_z]$:
$$
P_1 + p_1 = \bar{n}_x^2, 
\quad 
P_2 + p_2 = \bar{n}_y^2, 
\quad 
P_3 + p_3 = \bar{n}_z^2.
$$

In particular, we wish to analyze $\|\bar{\mathbf{n}}_{xy}\|^2 = \bar{n}_x^2 + \bar{n}_y^2$, which returns $(S_2^*, s_2^*)$, and also $\norm{\bar{\n}}^2$, which returns $(S_3^*, s_3^*)$ in \Cref{alg:sum_of_squares},

\begin{align*}
|S_2^*+s_2^*-\norm{\bar{\n}_{xy}}^2|&=|(P_1+p_1+P_2+p_2)(1+\varepsilon_{\text{SNN}}) - P_1-p_1-P_2-p_2|\\
&= |\varepsilon_{\text{SNN}}(P_1+p_1+P_2+p_2)|\\
\UB{|S_2^*+s_2^*-\left(\bar{n}_x^2+\bar{n}_y^2\right)|}&= 3u^2\norm{\bar{\n}_{xy}}^2
\end{align*}
\begin{align*}
|S_3^*+s_3^* - ||\mathbf{\bar{n}}||^2|&=|((P_1+p_1+P_2+p_2)(1+\varepsilon_{\text{SNN}})+P_3+p_3)(1+\varepsilon_{\text{SNN}})\\
&\quad\quad-(P_1+p_1+P_2+p_2+P_3+p_3)|\\
&=|(P_1+p_1)((1+\varepsilon_{\text{SNN}})^2-1)+\\
&\quad \; \;(P_2+p_2)((1+\varepsilon_{\text{SNN}})^2-1)+\\
&\quad \;\;(P_3+p_3)(1+\varepsilon_{\text{SNN}}-1)|\\
\UB{|S_3^*+s_3^* - ||\mathbf{\bar{n}}||^2|}&=\bar{n}_x^2(9u^4+6u^2)+ \bar{n}_y^2(9u^4+6u^2) + \bar{n}_z^2\cdot3u^2\\
&\leq 6u^2\norm{\bar{\n}}^2+\texttt{h.o.t}
\end{align*}

Building upon the previous analysis, we now examine \texttt{SumOfSquaresC}, which is used to compute $||\tilde{\n}_{xy}||^2 = \tilde{n}_x^2 + \tilde{n}_y^2$ and $||\tilde \n||^2 = \tilde{n}_x^2 + \tilde{n}_y^2 + \tilde{n}_z^2$  (including the error compensation term for the normal components). We use the subscripts 2 and 3 on the results ($S + s$) and intermediate values ($R^*$ and $R$) for the computation of $||\tilde \n_{xy}||^2$ and $||\tilde \n||^2$, respectively.
We also denote by $\varepsilon_{\text{CD}}$ the relative error associated with \Cref{alg:compensated_dot_product}.  
All subsequent bounds are derived using Wolfram Mathematica.

\paragraph{Upper Bound Derivations}

    \begin{align*}
        S_2 + s_2 &= S_2^* + R_2 \\
        &= S_2^* + \left(2R^*_2 + s_2^*\right)\left(1 + \delta_{14}\right) \\
        &= S_2^* + s_2^* + 2\left(\bar{n}_x e_{\bar{n}_x} + \bar{n}_y e_{\bar{n}_y}\right) 
        + 2\left(\bar{n}_x e_{\bar{n}_x} + \bar{n}_y e_{\bar{n}_y}\right)\varepsilon_{\text{CD}} \\
        &\quad
        + \left(2R^*_2 + s_2^*\right)\delta_{14} \\
        \UB{S_2 + s_2}&= \norm{\n_{xy}}^2+ 6u^2 \left(\norm{\n_{xy}}^2+2\norm{\n_{xy}}_1\right) + O\left(u^3\right)\\
        S_3 +s_3 &= S_3^* + R_3 \\
        &= S_3^* + \left(2R^*_3 + s_3^*\right)\left(1 + \delta_{15}\right) \\
       \UB{S_3 + s_3} &= \norm{\n}^2 + 3u^2 \left(3\norm{\n}^2+4\norm{\n}_1\right) 
        + O\left(u^3\right)
    \end{align*}
    \label{eq:S3_upperbound}

\paragraph{Lower Bound Derivations}

    \begin{align*}
        S_2 + s_2 &= S_2^* + R_2 \\
        \LB{S_2 + s_2}&=\norm{\n_{xy}}^2 - 2u^2 \left(3\norm{\n_{xy}}^2+2\norm{\n_{xy}}_1+4\right) + O\left(u^3\right)\\
        S_3 + s_3 &= S_3^* + R_3 \\
        \LB{S_3 + s_3}&=  \norm{\n}^2 - 3u^2 \left( 3\norm{\n}^2+4 \norm{\n}_1  + 4\right) 
        + O\left(u^3\right)
    \end{align*}
    \label{eq:S3_lowerbound}

%% file: main.bbl
\begin{thebibliography}{47}
\providecommand{\natexlab}[1]{#1}
\providecommand{\url}[1]{\texttt{#1}}
\expandafter\ifx\csname urlstyle\endcsname\relax
  \providecommand{\doi}[1]{doi: #1}\else
  \providecommand{\doi}{doi: \begingroup \urlstyle{rm}\Url}\fi

\bibitem[iee(1985)]{ieee754}
{IEEE} standard for binary {Floating-Point} arithmetic, 1985.

\bibitem[MPF(2023)]{MPFRWebsite}
{The MPFR Library}, 2023.
\newblock URL \url{https://www.mpfr.org}.
\newblock Accessed: June 30, 2023.

\bibitem[Berberich et~al.(2010)Berberich, Fogel, Halperin, Kerber, and
  Setter]{berberich2010arrangements}
Eric Berberich, Efi Fogel, Dan Halperin, Michael Kerber, and Ophir Setter.
\newblock Arrangements on parametric surfaces ii: Concretizations and
  applications.
\newblock \emph{Mathematics in Computer Science}, 4\penalty0 (1):\penalty0
  67--91, 2010.

\bibitem[Biscani and contributors(2024)]{mppp}
Francesco Biscani and contributors.
\newblock {mp++}: A high-performance multiprecision library, 2024.
\newblock URL \url{https://bluescarni.github.io/mppp/}.

\bibitem[Cherchi et~al.(2020)Cherchi, Livesu, Scateni, and Attene]{Cherchi2020}
Gianmarco Cherchi, Marco Livesu, Riccardo Scateni, and Marco Attene.
\newblock Fast and robust mesh arrangements using floating-point arithmetic.
\newblock \emph{ACM Trans. Graph.}, 39\penalty0 (6), November 2020.
\newblock ISSN 0730-0301.
\newblock \doi{10.1145/3414685.3417818}.
\newblock URL \url{https://doi.org/10.1145/3414685.3417818}.

\bibitem[de~Castro et~al.(2024)de~Castro, Cazals, Loriot, and
  Teillaud]{cgal:cclt-sgk3-24b}
Pedro Machado~Manh{\~a}es de~Castro, Fr{\'e}d{\'e}ric Cazals, S{\'e}bastien
  Loriot, and Monique Teillaud.
\newblock {3D} spherical geometry kernel.
\newblock In \emph{{CGAL} User and Reference Manual}. CGAL Editorial Board, 6.0
  edition, 2024.
\newblock URL
  \url{https://doc.cgal.org/6.0/Manual/packages.html#PkgCircularKernel3}.

\bibitem[Dekker(1971)]{dekker1971floating}
T.~J. Dekker.
\newblock A floating-point technique for extending the available precision.
\newblock \emph{Numerische Mathematik}, 18\penalty0 (3):\penalty0 224--242,
  1971.
\newblock \doi{10.1007/BF01397083}.

\bibitem[Dennis and Walker(1984)]{DennisWalker1984}
J.~E. Dennis and Homer~F. Walker.
\newblock Inaccuracy in quasi-{Newton} methods: Local improvement theorems.
\newblock In Bernhard Korte and Klaus Ritter, editors, \emph{Mathematical
  Programming at Oberwolfach II}, pages 70--85. Springer Berlin Heidelberg,
  Berlin, Heidelberg, 1984.
\newblock \doi{10.1007/BFb0121009}.

\bibitem[Devillers and Teillaud(2003)]{Devillers2003}
Olivier Devillers and Monique Teillaud.
\newblock Perturbations and vertex removal in a {3D} {Delaunay} triangulation.
\newblock In \emph{Proceedings of the Fourteenth Annual ACM-SIAM Symposium on
  Discrete Algorithms}, SODA '03, page 313–319, USA, 2003. Society for
  Industrial and Applied Mathematics.
\newblock ISBN 0898715385.

\bibitem[Edelsbrunner and M\"{u}cke(1990)]{Edelsbrunner1990}
Herbert Edelsbrunner and Ernst~Peter M\"{u}cke.
\newblock Simulation of {Simplicity}: a technique to cope with degenerate cases
  in geometric algorithms.
\newblock \emph{ACM Trans. Graph.}, 9\penalty0 (1):\penalty0 66–104, January
  1990.
\newblock ISSN 0730-0301.
\newblock \doi{10.1145/77635.77639}.
\newblock URL \url{https://doi.org/10.1145/77635.77639}.

\bibitem[Farrell et~al.(2009)Farrell, Piggott, Pain, Gorman, and
  Wilson]{FARRELL20092632}
P.E. Farrell, M.D. Piggott, C.C. Pain, G.J. Gorman, and C.R. Wilson.
\newblock Conservative interpolation between unstructured meshes via supermesh
  construction.
\newblock \emph{Computer Methods in Applied Mechanics and Engineering},
  198\penalty0 (33):\penalty0 2632--2642, 2009.
\newblock ISSN 0045-7825.
\newblock \doi{https://doi.org/10.1016/j.cma.2009.03.004}.
\newblock URL
  \url{https://www.sciencedirect.com/science/article/pii/S0045782509001315}.

\bibitem[Fousse et~al.(2007)Fousse, Hanrot, Lef\`{e}vre, P\'{e}lissier, and
  Zimmermann]{MPFR}
L.~Fousse, G.~Hanrot, V.~Lef\`{e}vre, P.~P\'{e}lissier, and P.~Zimmermann.
\newblock {MPFR}: A multiple-precision binary floating-point library with
  correct rounding, June 2007.

\bibitem[Graillat(2009)]{Graillat2009}
Stef Graillat.
\newblock Accurate floating-point product and exponentiation.
\newblock \emph{IEEE Transactions on Computers}, 58\penalty0 (7):\penalty0
  994--1000, 2009.
\newblock \doi{10.1109/TC.2008.215}.

\bibitem[Graillat et~al.(2015)Graillat, Lauter, Tang, Yamanaka, and
  Oishi]{Graillat2015}
Stef Graillat, Christoph Lauter, PING Tak~Peter Tang, Naoya Yamanaka, and
  Shin’ichi Oishi.
\newblock Efficient calculations of faithfully rounded l2-norms of n-vectors.
\newblock \emph{ACM Transactions on Mathematical Software}, 41\penalty0 (4),
  2015.
\newblock \doi{10.1145/2699469}.

\bibitem[Granlund and the GMP~development team(2023)]{GMP}
T.~Granlund and the GMP~development team.
\newblock {GNU MP}: The {GNU} multiple precision arithmetic library, 2023.
\newblock URL \url{http://gmplib.org/}.

\bibitem[Guennebaud et~al.(2010)Guennebaud, Jacob, et~al.]{eigenweb}
Ga\"{e}l Guennebaud, Beno\^{i}t Jacob, et~al.
\newblock Eigen v3.
\newblock http://eigen.tuxfamily.org, 2010.

\bibitem[Hanke et~al.(2016)Hanke, Redler, Holfeld, and Yastremsky]{YAC}
M.~Hanke, R.~Redler, T.~Holfeld, and M.~Yastremsky.
\newblock {YAC} 1.2.0: New aspects for coupling software in earth system
  modelling.
\newblock \emph{Geoscientific Model Development}, 9\penalty0 (8):\penalty0
  2755--2769, 2016.
\newblock \doi{10.5194/gmd-9-2755-2016}.
\newblock URL \url{https://gmd.copernicus.org/articles/9/2755/2016/}.

\bibitem[Higham(2002)]{Higham2002}
Nicholas~J. Higham.
\newblock \emph{Accuracy and Stability of Numerical Algorithms}.
\newblock Society for Industrial and Applied Mathematics, 2nd edition, 2002.
\newblock \doi{10.1137/1.9780898718027}.
\newblock URL \url{https://epubs.siam.org/doi/abs/10.1137/1.9780898718027}.

\bibitem[Hu et~al.(2018)Hu, Zhou, Gao, Jacobson, Zorin, and
  Panozzo]{Hu2018TetMeshWild}
Yixin Hu, Qingnan Zhou, Xifeng Gao, Alec Jacobson, Denis Zorin, and Daniele
  Panozzo.
\newblock Tetrahedral meshing in the wild.
\newblock \emph{ACM Trans. Graph.}, 37\penalty0 (4), July 2018.
\newblock ISSN 0730-0301.
\newblock \doi{10.1145/3197517.3201353}.
\newblock URL \url{https://doi.org/10.1145/3197517.3201353}.

\bibitem[Hu et~al.(2020)Hu, Schneider, Wang, Zorin, and Panozzo]{hu_fast_2020}
Yixin Hu, Teseo Schneider, Bolun Wang, Denis Zorin, and Daniele Panozzo.
\newblock Fast tetrahedral meshing in the wild.
\newblock \emph{ACM Transactions on Graphics}, 39\penalty0 (4), August 2020.
\newblock ISSN 0730-0301, 1557-7368.
\newblock \doi{10.1145/3386569.3392385}.
\newblock URL \url{https://dl.acm.org/doi/10.1145/3386569.3392385}.

\bibitem[{Intel and UXL Foundation}(2025)]{intel_onetbb_2022_2_0}
{Intel and UXL Foundation}.
\newblock {Intel {oneAPI} Threading Building Blocks {(oneTBB)}}, 2025.
\newblock URL \url{https://github.com/oneapi-src/oneTBB}.

\bibitem[Jeannerod and Rump(2017)]{JeannerodRump2017}
Claude-Pierre Jeannerod and Siegfried~M. Rump.
\newblock On relative errors of floating-point operations: Optimal bounds and
  applications.
\newblock \emph{Mathematics of Computation}, 87:\penalty0 803--819, 2017.
\newblock URL \url{https://api.semanticscholar.org/CorpusID:40376156}.

\bibitem[Jeannerod et~al.(2013)Jeannerod, Louvet, and Muller]{KahanAlgo}
Claude-Pierre Jeannerod, Nicolas Louvet, and Jean-Michel Muller.
\newblock Further analysis of {Kahan's} algorithm for the accurate computation
  of 2 x 2 determinants.
\newblock \emph{Mathematics of Computation}, 82\penalty0 (284):\penalty0
  2245--2264, 2013.
\newblock \doi{10.1090/S0025-5718-2013-02679-8}.
\newblock URL \url{https://ens-lyon.hal.science/ensl-00649347v4}.

\bibitem[Kahan(1996)]{KahanLectureNote}
W.~Kahan.
\newblock Lecture notes on the status of {IEEE-754}, 1996.
\newblock URL
  \url{http://www.cs.berkeley.edu/~wkahan/ieee754status/IEEE754.PDF}.

\bibitem[Kahan(2004)]{kahancost}
W~Kahan.
\newblock On the cost of floating-point computation without extra-precise
  arithmetic. 2004.
\newblock \emph{URL: http://www. cs. berkeley. edu/-wkahan/Qdrtcs. pdf}, 2004.

\bibitem[Knuth(1997)]{Knuth1997}
Donald~E. Knuth.
\newblock \emph{The {Art of Computer Programming}, Volume 2 (3rd Ed.):
  Seminumerical Algorithms}.
\newblock Addison-Wesley Longman Publishing Co., Inc., USA, 1997.
\newblock ISBN 0201896842.

\bibitem[Krumm(2000)]{krumm2000}
John Krumm.
\newblock Intersection of two planes.
\newblock May 2000.
\newblock URL
  \url{https://www.microsoft.com/en-us/research/publication/intersection-of-two-planes/}.

\bibitem[Lef\`{e}vre et~al.(2023)Lef\`{e}vre, Louvet, Muller, Picot, and
  Rideau]{LefevreEtAl2023}
Vincent Lef\`{e}vre, Nicolas Louvet, Jean-Michel Muller, Joris Picot, and
  Laurence Rideau.
\newblock Accurate calculation of {Euclidean} norms using double-word
  arithmetic.
\newblock \emph{ACM Transactions on Mathematical Software}, 49\penalty0
  (1):\penalty0 1--34, March 2023.
\newblock \doi{10.1145/3568672}.

\bibitem[Marsico and Ullrich(2023)]{strategiesRemapping}
D.~H. Marsico and P.~A. Ullrich.
\newblock Strategies for conservative and non-conservative monotone remapping
  on the sphere.
\newblock \emph{Geoscientific Model Development}, 16\penalty0 (5):\penalty0
  1537--1551, 2023.
\newblock \doi{10.5194/gmd-16-1537-2023}.
\newblock URL \url{https://gmd.copernicus.org/articles/16/1537/2023/}.

\bibitem[Mirshak(2024)]{mirshak2024intersections}
Ramzi Mirshak.
\newblock On finding the intersections of a small circle with a great circle
  segment on a sphere: Methods and code.
\newblock Reference Document DRDC-RDDC-2024-D041, Defence Research and
  Development Canada (DRDC), Ottawa, Canada, May 2024.
\newblock URL
  \url{https://cradpdf.drdc-rddc.gc.ca/PDFS/unc471/p818157_A1b.pdf}.
\newblock This document has been reviewed and does not contain controlled
  technical data. Approved for public release.

\bibitem[Muller et~al.(2018)Muller, Brunie, de~Dinechin, Jeannerod, Joldes,
  Lef\`{e}vre, Melquiond, Revol, and Torres]{MullerEtAl2018}
Jean-Michel Muller, Nicolas Brunie, Florent de~Dinechin, Claude-Pierre
  Jeannerod, Mioara Joldes, Vincent Lef\`{e}vre, Guillaume Melquiond, Nathalie
  Revol, and Serge Torres.
\newblock \emph{Handbook of Floating-Point Arithmetic}.
\newblock Birkh\"auser Boston, 2nd edition, 2018.
\newblock ISBN 978-3-319-76525-9.

\bibitem[Ogita et~al.(2005)Ogita, Rump, and Oishi]{Ogita2005}
Takeshi Ogita, Siegfried~M. Rump, and Shin'ichi Oishi.
\newblock Accurate sum and dot product.
\newblock \emph{SIAM Journal on Scientific Computing}, 26\penalty0
  (6):\penalty0 1955--1988, 2005.
\newblock \doi{10.1137/030601818}.

\bibitem[Project(2024)]{libquadmath}
GNU Project.
\newblock {GNU} libquadmath manual, 2024.
\newblock URL \url{https://gcc.gnu.org/onlinedocs/libquadmath/}.
\newblock Accessed: 2024-07-18.

\bibitem[Richard~Shewchuk(1997)]{richard_shewchuk_adaptive_1997}
Jonathan Richard~Shewchuk.
\newblock Adaptive precision floating-point arithmetic and fast robust
  geometric predicates.
\newblock \emph{Discrete \& Computational Geometry}, 18\penalty0 (3):\penalty0
  305--363, October 1997.
\newblock ISSN 0179-5376.
\newblock \doi{10.1007/PL00009321}.
\newblock URL \url{http://link.springer.com/10.1007/PL00009321}.

\bibitem[Rump(2023)]{Rump2023}
Siegfried Rump.
\newblock Fast and accurate computation of the {Euclidean} norm of a vector.
\newblock \emph{Japan Journal of Industrial and Applied Mathematics}, 40, 2023.
\newblock \doi{10.1007/s13160-023-00593-8}.

\bibitem[Rump(2009)]{rump2009error}
Siegfried~M Rump.
\newblock Error-free transformations and ill-conditioned problems.
\newblock In \emph{International Workshop on Verified Computations and Related
  Topics}, 2009.

\bibitem[Schling(2011)]{boost}
Boris Schling.
\newblock \emph{The {Boost C++} Libraries}.
\newblock XML Press, 2011.
\newblock ISBN 0982219199.

\bibitem[Schulzweida(2023)]{CDO}
Uwe Schulzweida.
\newblock {CDO User Guide}, October 2023.
\newblock URL \url{https://doi.org/10.5281/zenodo.10020800}.

\bibitem[Shewchuk(1996)]{Triangle96}
Jonathan~Richard Shewchuk.
\newblock {Triangle}: Engineering a {2D} quality mesh generator and delaunay
  triangulator.
\newblock In \emph{{WACG}}, volume 1148 of \emph{Lecture Notes in Computer
  Science}, pages 203--222. Springer, 1996.

\bibitem[{The CGAL Project}(2024)]{cgal:eb-24b}
{The CGAL Project}.
\newblock \emph{{CGAL} User and Reference Manual}.
\newblock CGAL Editorial Board, 6.0 edition, 2024.
\newblock URL \url{https://doc.cgal.org/6.0/Manual/packages.html}.

\bibitem[Trettner et~al.(2022)Trettner, Nehring-Wirxel, and
  Kobbelt]{trettner_ember_2022}
Philip Trettner, Julius Nehring-Wirxel, and Leif Kobbelt.
\newblock {EMBER}: Exact mesh booleans via efficient \& robust local
  arrangements.
\newblock \emph{ACM Transactions on Graphics}, 41\penalty0 (4):\penalty0 1--15,
  July 2022.
\newblock ISSN 0730-0301, 1557-7368.
\newblock \doi{10.1145/3528223.3530181}.
\newblock URL \url{https://dl.acm.org/doi/10.1145/3528223.3530181}.

\bibitem[Ullrich and Taylor(2015)]{TempestRemap}
Paul~A. Ullrich and Mark~A. Taylor.
\newblock Arbitrary-order conservative and consistent remapping and a theory of
  linear maps: Part i.
\newblock \emph{Monthly Weather Review}, 143\penalty0 (6):\penalty0 2419 --
  2440, 2015.
\newblock \doi{10.1175/MWR-D-14-00343.1}.
\newblock URL
  \url{https://journals.ametsoc.org/view/journals/mwre/143/6/mwr-d-14-00343.1.xml}.

\bibitem[Ullrich et~al.(2016)Ullrich, Devendran, and Johansen]{TempestRemap2}
Paul~A. Ullrich, Dharshi Devendran, and Hans Johansen.
\newblock Arbitrary-order conservative and consistent remapping and a theory of
  linear maps: Part ii.
\newblock \emph{Monthly Weather Review}, 144\penalty0 (4):\penalty0 1529 --
  1549, 2016.
\newblock \doi{10.1175/MWR-D-15-0301.1}.
\newblock URL
  \url{https://journals.ametsoc.org/view/journals/mwre/144/4/mwr-d-15-0301.1.xml}.

\bibitem[{UXarray Organization}(2024)]{uxarray2024}
{UXarray Organization}.
\newblock {UXarray} (v2024.03.0) [software], 2024.
\newblock URL \url{https://doi.org/10.5281/zenodo.10895493}.
\newblock Project Raijin \& Project SEATS.

\bibitem[Wein et~al.(2024)Wein, Berberich, Fogel, Halperin, Hemmer, Salzman,
  and Zukerman]{cgal:wfzh-a2-24b}
Ron Wein, Eric Berberich, Efi Fogel, Dan Halperin, Michael Hemmer, Oren
  Salzman, and Baruch Zukerman.
\newblock {2D} arrangements.
\newblock In \emph{{CGAL} User and Reference Manual}. {CGAL Editorial Board},
  {6.0} edition, 2024.
\newblock URL
  \url{https://doc.cgal.org/6.0/Manual/packages.html#PkgArrangementOnSurface2}.

\bibitem[{Wolfram Research, Inc.}(2024)]{Mathematica}
{Wolfram Research, Inc.}
\newblock {Mathematica}, version 14.0, 2024.
\newblock URL \url{https://www.wolfram.com/mathematica}.

\bibitem[Zhou et~al.(2016)Zhou, Grinspun, Zorin, and Jacobson]{Zhou:2016:MASG}
Qingnan Zhou, Eitan Grinspun, Denis Zorin, and Alec Jacobson.
\newblock Mesh arrangements for solid geometry.
\newblock \emph{ACM Transactions on Graphics (TOG)}, 35\penalty0 (4), 2016.

\end{thebibliography}
